\documentclass[11pt]{article}
\catcode`\@=11
\catcode`\@=12

\usepackage{amsmath}
\usepackage{amssymb}
\usepackage{amsthm}

\newcommand{\Rm}{\mathbb{R}}

\newcommand{\NN}{\mathbb{N}}

\newcommand{\mC}{\ensuremath{\mathcal{C}}}

\newcommand{\mT}{\ensuremath{\mathcal{T}}}

\newcommand{\mF}{\ensuremath{\mathcal{F}}}
\newcommand{\Nm}{\ensuremath{\mathbb{N}}}

\newcommand{\Hm}{\ensuremath{\mathbb{H}}}
\newcommand{\mM}{\ensuremath{\mathcal{M}}}
\newcommand{\mK}{\ensuremath{\mathcal{K}}}

\newcommand{\mO}{\ensuremath{\mathcal{O}}}
\newcommand{\mI}{\ensuremath{\mathcal{I}}}
\newcommand{\mB}{\ensuremath{\mathcal{B}}}

\newcommand{\mE}{\ensuremath{\mathcal{E}}}
\newcommand{\vs}{\vspace{.5cm}}

\newtheorem{lem}{Lemma}
\newtheorem{thm}{Theorem}

\newtheorem{prop}[lem]{Proposition}
\newtheorem{defn}[lem]{Definition}

\def\proof {\noindent{\sc{Proof. }}}
\def\qed {\mbox{}\hfill {\small \fbox{}} \\}  
\def\lto{\longrightarrow}
\def\lmto{\longmapsto}

\def\leq{\leqslant}
\def\geq{\geqslant}

\title{Optimal mass transportation and Mather theory}
\author{Patrick  Bernard and Boris  Buffoni}
\date{July 2005 }
\begin{document}

\maketitle

Abstract:
We study the Monge transportation problem when
the cost is the action associated to a Lagrangian function
on a compact manifold.
We show that the transportation can be interpolated by
a Lipschitz lamination. We describe several direct
variational problems the minimizers of which are these
Lipschitz laminations.
We prove the existence of an  optimal transport map
when the transported measure is absolutely continuous.
We explain the relations with Mather's minimal measures.
\\
\vs

R\'esum\'e:
On \'etudie le probl\`eme de transport de Monge lorsque 
le cout est l'action associ\'ee \`a un Lagrangien sur une vari\'et\'e
compacte.
On montre que le transport peut \^etre interpol\'e par
une lamination lipschitzienne.
On d\'ecrit plusieurs probl\`emes variationnels 
directs dont ces laminations sont les minimiseurs.
On montre l'existence d'une application de transport 
optimale lorsque la mesure transport\'ee est absolument continue.
On explique les relations avec les mesures minimisantes de Mather.
\\
\vs
\\

{\rm Patrick Bernard}\\
Institut Fourier, Grenoble,\\
on move to \\
CEREMADE\\
Universit\'e de Paris Dauphine\\
Pl. du Mar\'echal de Lattre de Tassigny\\
75775 Paris Cedex 16\\
France\\
\texttt{patrick.bernard@ceremade.dauphine.fr}\\

\vs

{\rm Boris Buffoni}\\
School of Mathematics\\
\'Ecole Polytechnique F\'ed\'erale-Lausanne\\
SB/IACS/ANA Station 8\\
1015 Lausanne\\
Switzerland\\
\texttt{boris.buffoni@epfl.ch}

\newpage

Several   observations have recently renewed the 
interest for the classical topic of optimal mass transportation,
whose primary origin is attributed to Monge
a few years before French revolution.
The framework is as follows. A space $M$ is given,
which in the present paper will be a compact manifold,
as well as a continuous cost function 
$c(x,y):M\times M \lto \Rm$.
Given two probability  measures $\mu_0$ and $\mu_1$ on $M$,
the mappings $\Psi:M\lto M$ which 
transport $\mu_0$ into $\mu_1$ and minimize
the total cost 
$\int_M c(x,\Psi(x))d\mu_0$
are studied.
It turns out, and it was the core of the investigations of Monge, 
 that these mappings have very remarkable geometric properties,
at least at a formal level.

Only much more recently was the question of the existence
of optimal objects rigorously solved  by Kantorovich 
in a famous paper of 1942. Here we speak of optimal objects,
and not of optimal mappings, because 
the question of existence of an optimal mapping is ill-posed, 
so that the notion of optimal objects has to be relaxed,
in a way that nowadays seems very natural, and that 
was discovered by  Kantorovich.

Our purpose here is to continue the work initiated
by Monge, recently awakened  by Brenier
and enriched by  other authors, on the study 
of geometric properties of optimal objects.
The costs functions we consider are natural generalizations 
of the cost $c(x,y)=d(x,y)^2$
considered by Brenier and many other authors.
The book \cite{Vi:03} gives some ideas of the 
applications expected from this kind of questions.
More precisely, we consider 
a Lagrangian function 
$L(x,v,t):TM\times \Rm\lto \Rm$
which is convex in $v$ and satisfies standard 
hypotheses recalled later, 
and define our cost by
$$
c(x,y)=\min_{\gamma} \int _0^1 L(\gamma(t),\dot \gamma(t),t) dt
$$
where the minimum is taken on the set of curves
$\gamma:[0,1]\lto M$ satisfying $\gamma(0)=x$
and $\gamma(1)=y$. 
Note that this class of costs does not contain the very natural
cost  $c(x,y)=d(x,y)$.
Such costs are studied in a  second paper \cite{BeBu:second}.

Our main result is that the optimal transports can be interpolated
by measured Lipschitz laminations, or geometric currents
in the sense of Ruelle and Sullivan.
Interpolations of transport have  already been considered by
Benamou,
Brenier and McCann for less general cost functions,
and with different purposes.
Our methods are inspired by the theory 
of Mather, Ma\~n\'e and Fathi on Lagrangian dynamics,
and we  will  detail rigorously the relations between
these theories. Roughly, they are exactly similar
except that mass transportation is a Dirichlet boundary value
problem, while Mather theory is a periodic boundary value
Problem. 
We will also prove,
extending works of Brenier, Gangbo, 
McCann, Carlier, and other authors,  
that the optimal transportation can be performed
by a Borel map with the additional assumption that the transported 
measure is absolutely continuous.

Various connections between Mather-Fathi theory, optimal
mass transportation and Hamilton-Jacobi equations
have recently been discussed, mainly at a formal level,
in the literature, see for example  \cite{Vi:03}, or \cite{EvGo:02}, where
they are all presented  as infinite dimensional linear 
programming problems.
This have motivated a lot of activity around the interface between 
Aubry-Mather theory and optimal transportation,
some of which  overlap partly the present work.
For example, 
at the moment of submitting the paper,
we have been informed of the existence 
of the recent preprints of De Pascale, Gelli and Granieri, 
\cite{PaStGr}, and of Granieri, \cite{Gr}.
We had also been aware of a manuscript  by  Wolansky \cite{Wo}
for a few weeks, which, independently,  
and by somewhat different methods, obtains results
similar to ours. 
Note however that 
Lipschitz regularity, which we consider as one of our 
most important results, was not obtained
in this preliminary version of \cite{Wo}.
It is worth also mentioning the papers 
 \cite{Pr:04}
of Pratelli and \cite{Lo:05} of Loeper.

This paper emanates from the collaboration of the Authors 
during the end  of the  stay of the first author  in EPFL for the 
academic year 2002-2003, granted by  
the Swiss National Science Foundation.

\section{Introduction}
We present the  context and the main results of
the paper.

\subsection{Lagrangian, Hamiltonian  and cost}
In all the present paper, the space $M$ will be
 a compact  and connected Riemannian manifold
without boundary. Some standing notations are gathered in the 
appendix.
Let us fix a positive real number $T$, and a Lagrangian function
$$
L\in C^2(TM\times [0,T],\Rm).
$$
A curve $\gamma\in C^2([0,T],M)$ is called an extremal
if it is a critical point of the action 
$$
\int_0^T L(\gamma(t),\dot \gamma(t),t)dt$$
 with fixed endpoints.
It is called a minimizing extremal if it is minimizing the action.
We assume:\\
\textbf{convexity}
For each $(x,t)\in M\times[0,T]$, 
the function $v\lmto L(x,v,t)$ is convex with positive definite Hessian
at each point.\\
\textbf{superlinearity}
For each  $(x,t)\in M\times[0,T]$, we have $L(x,v,t)/\|v\|\lto \infty$
as $\|v\|\lto \infty$. Arguing as in \cite{Fa:un} (Lemma 3.2.2), 
this implies that
for all $\alpha>0$ there exists $C>0$ such that $L(x,v,t)\geq \alpha \|v\|-C$
for all $(x,v,t)\in TM\times[0,T]$.\\
\textbf{completeness}
For each $(x,v,t)\in TM \times[0,T]$, there exists one and only one 
extremal $\gamma\in C^2([0,T],M)$ such that 
$(\gamma(t),\dot \gamma(t))=(x,v)$.

We associate to the Lagrangian $L$ a Hamiltonian function
$H\in C^2(T^*M\times[0,T],\Rm)$ given by
$$
H(x,p,t)=\max_v p(v)-L(x,v,t).
$$
We endow the cotangent bundle   $T^*M$  with its canonical symplectic structure,
and  associate to the Hamiltonian $H$
the time-dependent vectorfield $Y$ on $T^*M$,
which is given by 
$$
Y=(\partial_p H, -\partial_x H)
$$
in any canonical local trivialisation of $T^*M$.
The hypotheses on $L$ can be expressed in terms of the function $H$:\\
\textbf{convexity}
For each $(x,t)\in M\times[0,T]$, 
the function $p\lmto H(x,p,t)$ is convex with positive definite Hessian
at each point.\\
\textbf{superlinearity}
For each  $(x,t)\in M\times[0,T]$, we have $H(x,p,t)/\|p\|\lto \infty$
as $\|p\|\lto \infty$.\\
\textbf{completeness}
Each solution of the equation $(\dot x(t),\dot p(t))=Y(x(t),p(t),t)$
can be extended to the interval $[0,T]$.
We can then define, 
for all $s,t\in[0,T]$, 
the flow $\varphi_s^t$ of $Y$ from times $s$ to time $t$.

In addition, the mapping
$\partial_v L:TM\times[0,T]\lto T^*M\times[0,T]$
is a $C^1$ diffeomorphism, whose inverse is the mapping
$\partial_p H$.
These diffeomorphisms  conjugate  $Y$ with a time-dependent
vectorfield $E$ on $TM$.
We denote the flow of $E$ by
$\psi_s ^t:TM\lto TM$ ($s,t \in [0,T]$), which is such that 
$\psi_s^s=Id$ and
$\partial_t \psi_s^t=E_t\circ \psi_s^t$,
where as usual $E_t$ denotes the vectorfield $E(.,t)$ on $TM$. 
The diffeomorphisms $\partial_v L$ and $\partial_p H$ conjugate the flows
$\psi_s^t$ and $\varphi_s^t$. Moreover the extremals are the projection 
of the integral curves of $E$ and
\begin{equation}\label{elflow}
\Big(\pi \circ \psi_s^t\,,\,\partial_t(\pi \circ \psi_s^t)\Big)=\psi_s^t,
\end{equation}
where $\pi:TM\rightarrow M$ is the canonical projection.
In \eqref{elflow}, $\partial_t(\pi \circ \psi_s^t)$ is seen
as a vector in the tangent space of $M$ at
$\pi \circ \psi_s^t$.  If
$\partial_t(\pi \circ \psi_s^t)$ 
is seen as a point in $TM$,
\eqref{elflow} becomes simply
$\partial_t(\pi \circ \psi_s^t)=\psi_s^t$.

For each $0\leq s <t\leq T$, we define the cost function 
$$
c_s^t(x,y)=\min_{\gamma} \int _s^t L(\gamma(\sigma),\dot \gamma(\sigma),
\sigma) d\sigma
$$
where the minimum is taken on the set of curves
$\gamma\in C^2([s,t],M)$ satisfying $\gamma(s)=x$
and $\gamma(t)=y$. 
That this minimum exists is a standard result
under our  hypotheses, see \cite{Ma:91} or \cite{Fa:un}.

\begin{prop}\label{extremals}
Let us fix a subinterval $[s,t]\subset [0,T]$.
The set $\mE\subset  C^2([s,t],M)$ of minimizing 
extremals is compact for the $C^2$ topology. 
\end{prop}

Let us mention that,
for each $(x_0,s)\in M\times[0,T]$, the function
$(x,t)\lmto c_s^t(x_0,x)$
is a viscosity solution of the Hamilton-Jacobi equation
$$
\partial_t u +
H(x,\partial_x u,t)=0
$$
on $M\times]s,T[$.
This remark may help the reader in understanding the key 
role which will be played by this equation in the sequel.


\subsection{Monge-Kantorovich theory}\label{MKT}
We recall the basics of Monge-Kantorovich duality.
The proofs are available in many texts on the subjects, for
example \cite{Am:00,RaRu:98,Vi:03}.
We assume that $M$ is a compact manifold and that 
$c(x,y)$ is a continuous cost function on 
$M\times M$, which will later be one of 
the costs $c_s^t$ defined above.
Given two Borel probability measures $\mu_0$
and $\mu_1$ on $M$,
a transport plan between $\mu_0$ and 
$\mu_1$ is a measure on $M\times M$
which satisfies 
$$\pi_{0\sharp}(\eta)=\mu_0 \text{ and }
\pi_{1\sharp}(\eta)=\mu_1,
$$
where $\pi_0:M\times M\lto M $ is the projection
on the first factor, and $\pi_1$ is the projection on
the second factor. 
We denote by $\mK(\mu_0,\mu_1)$, after Kantorovich,
the set of transport plans.
Kantorovich proved the existence of 
a minimum in the expression 
$$
C(\mu_0,\mu_1)=
\min_{\eta\in \mK(\mu_0,\mu_1)} \int _{M\times M} c d\eta
$$
for each  pair $(\mu_0, \mu_1)$ of probability measures on $M$.
Here we will denote by 
\begin{equation}\label{Monge}
C_s^t(\mu_0,\mu_1)
:= \min _{\eta\in \mK(\mu_0,\mu_1)}  \int_{M\times M}
c_s^t(x,y)d\eta(x,y)
\end{equation}
the optimal value associated to our family of costs $c_s^t$
The plans which realize this minimum are called
optimal transfer plans.
A pair $(\phi_0,\phi_1)$ of continuous functions 
is called an admissible Kantorovich pair if is satisfies the 
relations
$$
\phi_1(x)=\min_{y\in M} \phi_0(y)+c(y,x)
\text{ and } 
\phi_0(x)=\max _{y\in M} \phi_1(y)-c(x,y)
$$
for all point  $x\in M$.
Note that the admissible pairs are composed of Lipschitz functions
if the cost $c$ is Lipschitz, 
which is the case of the costs
$c_s^t$ when $s<t$.
Another discovery of Kantorovich is that 
\begin{equation}\label{Kanto}
C(\mu_0,\mu_1)
=\max _{\phi_0,\phi_1}
\Big(
\int_M \phi_1d\mu_1-\int_M \phi_0d\mu_0
\Big)
\end{equation}
where the maximum is taken on the set of 
 admissible Kantorovich pairs $(\phi_0, \phi_1)$.
This maximization problem is called the dual Kantorovich problem,
the admissible pairs  which reach this maximum
are called optimal Kantorovich pairs.
The direct problem (\ref{Monge}) and dual problem (\ref{Kanto})
are related as follows. 
\begin{prop}
If $\eta$ is an optimal transfer plan, and if $(\phi_0,\phi_1)$
is a Kantorovich  optimal pair, then the support of $\eta$ 
is contained in the set 
$$
\{(x,y)\in M^2 \text{ such that }
\phi_1(y)-\phi_0(x)=c(x,y)\}\subset M\times M.
$$
\end{prop}
Let us remark that the knowledge of the set of 
Kantorovich admissible pairs
is equivalent to the knowledge of the cost function $c$.
\begin{lem}
We have 
$$c(x,y)=\max _{(\phi_0,\phi_1)} \phi_1(y)-\phi_0(x)
$$
where the maximum is taken on the set of Kantorovich
admissible pairs.
\end{lem}
\proof
This maximum is clearly less that $c(x,y)$.
For the other inequality, let us 
fix  points $x_0$ and $y_0$ in $M$, and consider the functions
$\phi_1(y)=c(x_0,y)$
and 
$\phi_0(x)=\max  _{y\in M}\phi_1(y) -c(x,y)$.
We have 
$\phi_1(y_0)-\phi_0(x_0)=
c(x_0,y_0)-0=
c(x_0,y_0)$.
So it is enough to prove that the pair $(\phi_0,\phi_1)$
is an admissible Kantorovich pair,
and more precisely that 
$\phi_1(y)=\min_{x\in M} \phi_0(x)+c(x,y)$.
We have 
$$\phi_0(x)+c(x,y)\geq c(x_0,y)-c(x,y)+c(x,y)\geq c(x_0,y)=
\phi_1(y)$$
which gives the inequality 
$\phi_1(y)\leq  \min_{x\in M} \phi_0(x)+c(x,y)$.
On the other hand, we have 
$$\min_{x\in M} \phi_0(x)+c(x,y)
\leq \phi_0(x_0)+c(x_0,y)=c(x_0,y)=\phi_1(y).
$$
\qed
\subsection{Interpolations}\label{subsection Interpolations}
In this section, the Lagrangian $L$ and time $T>0$
are fixed.
It is not hard to see that, if $\mu_1,\mu_2$ and $\mu_3$
are three probability measures on $M$, and if 
$t_1\leq t_2 \leq t_3\in [0,T]$ are three times, then we have the inequality
$$
C_{t_1}^{t_3}(\mu_{1},\mu_{3})\leq  
C_{t_1}^{t_2}(\mu_{1},\mu_{2})+ 
C_{t_2}^{t_3}(\mu_{2},\mu_{3}). 
$$
The family $\mu_t,t\in [0,T]$ of probability measures
on $M$ is called an interpolation between $\mu_0$ and
$\mu_T$ if it satisfies the equality
$$
C_{t_1}^{t_3}(\mu_{t_1},\mu_{t_3})= 
C_{t_1}^{t_2}(\mu_{t_1},\mu_{t_2})+ 
C_{t_2}^{t_3}(\mu_{t_2},\mu_{t_3}) 
$$
for all $0 \leq t_1\leq t_2 \leq t_3\leq T$.
Our main result is the following:

\noindent
\textbf{Theorem A.}
\begin{itshape}
For each pair $\mu_0, \mu_T$
of probability measures, there exist interpolations
between $\mu_0$ and $\mu_T$.
Moreover, each interpolation $(\mu_t),t\in [0,T]$
is given by a Lipschitz measured lamination in the following sense:

\noindent
\textbf{\emph{Eulerian description :}}
There exists a bounded  locally Lipschitz vectorfield 
$X(x,t):M\times ]0,T[\lto TM$ such that,
if $\Psi_s^t,(s,t)\in ]0,T[^2$ is the flow of $X$
from time $s$ to time $t$, then
$(\Psi_s^t)_{\sharp}\mu_s=\mu_t$
for each $(s,t)\in ]0,T[^2$.

\noindent
\textbf{\emph{Lagrangian description :}}
There exists a family $\mF\subset  C^2([0,T],M)$ of minimizing extremals
$\gamma$ of $L$, which is such that 
the relation $\dot \gamma(t)=X(\gamma(t),t)$ holds
for each $t\in ]0,T[$ and for each $\gamma \in \mF$.
The set 
$$\tilde \mT=\{(\gamma(t),\dot \gamma(t),t), t\in ]0,T[,\gamma\in \mF\}
\subset TM\times]0,T[
$$
is invariant under the Euler-Lagrange flow $\psi$.
The measure $\mu_t$ is supported on 
$
\mT_t=\{\gamma(t),\gamma\in \mF\}.
$
In addition, 
there exists a   continuous family $m_t,t\in[0,T]$ of probability
measures on $TM$ such that $m_t$ is concentrated 
on 
$
\tilde \mT_t=\{(\gamma(t),\dot \gamma(t)),\gamma\in \mF\}
$ 
for each $t\in ]0,T[$, such that
$\pi_{\sharp} m_t=\mu_t$ for each $t\in [0,T]$, and such that 
$$
m_t=(\psi_s^t)_{\sharp}m_s
$$
for all $(s,t)\in [0,T]^2$.

\noindent
\textbf{\emph{Hamilton-Jacobi equation :}}
There exists a Lipschitz and  $C^1$ function $v(x,t):M\times]0,T[\lto \Rm$
 which satisfies  the inequation 
$$
\partial_tv+H(x,\partial_xv,t)\leq 0,
$$
with equality if and only if 
$(x,t)\in \mT=\{(\gamma(t),t),\gamma\in \mF, t\in ]0,T[\}$,
and such that 
$X(x,t)=\partial_pH(x,\partial_x v(x,t),t)$
for each $(x,t)\in \mT$.

\noindent
\textbf{\emph{Uniqueness  :}}
There may exist several different interpolations.
However, it is possible to choose the vectorfield $X$, the 
family $\mF$ and the sub-solution $v$ in such a way   that the statements 
above hold for all interpolations $\mu_t$ with these fixed
$X$, $\mF$ and $v$.
For each $s<t\in ]0,T[$,
the measure $(Id\times \Psi_s^t)_{\sharp} \mu_s$
is the only optimal transport plan in $\mK(\mu_s,\mu_t)$
for the cost $c_s^t$. This implies that 
$$
\int _M c_s^t(x,\Psi_s^t(x)) d\mu_s(x)=C_s^t(\mu_s,\mu_t).
$$
\end{itshape}

Let us comment a bit the preceding statement.
The set $\tilde \mT\subset TM\times]0,T[$ is the image by 
the Lipschitz map $(x,t)\rightarrow (X(x,t),t)$
of the set  $\mT\subset TM\times]0,T[$.
We shall not take 
$X(x,t)=\partial_pH(x,\partial_xv(x,t),t)$
outside of $\mT$ because we do not prove that 
this vectorfield is Lipschitz outside of $\mT$.
The data of the vectorfield $X$ outside of $\mT$
is immaterial: any Lipschitz extension of $X_{|\mT}$
will fit. Note also that the relation 
\begin{equation}
\label{X and flow}
\Psi_s^t=\pi\circ \psi_s^t \circ X_s
\end{equation}
holds on $\mT_s$, where $X_s(.)=X(.,s)$.

The vectorfield $X$ in the statement depends on 
the transported measures $\mu_0$ and $\mu_T$.
The Lipschitz constant of $X$, however, can be fixed
independently of these measures, as we now state
(see Proposition \ref{Fathilip}, 
Proposition \ref{mT},
Theorem \ref{construction of v}
and \eqref{refX}):\vs

\noindent
\textbf{Addendum}
\begin{itshape}
There exists a decreasing function $K(\epsilon):]0,T/2[\lto ]0,\infty[$,
which depends only on the time $T$ and on the Lagrangian $L$, and
such that, for each pair $\mu_0, \mu_T$ of probability measures,
one can choose the vectorfield $X$ in Theorem A in such a way that 
$X$ is $K(\epsilon)$-Lipschitz on 
$M\times [\epsilon,T-\epsilon]$
for each $\epsilon\in ]0,T/2[$.
\end{itshape}\vs

Proving Theorem A is the main goal of the present paper.
We will present in section \ref{problems} some direct variational 
problems which are well-posed and of which  the transport 
interpolations are in some sense the solutions.
We believe that these variational problems are interesting
in themselves.
In order to describe the solutions of the variational problem,
we will rely on a dual approach based on the Hamilton-Jacobi equation,
inspired from Fathi's approach to Mather theory,
as detailed in section \ref{sectionHJ}.
The solutions of the problems of section \ref{problems},
as well as the transport interpolations, are then described
in section \ref{ODP}, which ends the proof of Theorem A.

\subsection{Case of an absolutely continuous measure $\mu_0$}
Additional conclusions concerning
optimal transport can usually be obtained when the initial measure 
$\mu_0$ is absolutely continuous.
For example a standard question is whether the optimal transport
can be realized by an optimal mapping.

A transport map is a Borel map $\Psi:M\lto M$ which satisfies
$\Psi _{\sharp}\mu_0=\mu_1$.
To any transport map $\Psi$ is naturally associated the 
transport plan 
$(Id \times \Psi)_{\sharp} \mu_0$, called the induced transport plan.
An optimal map is a transport map $\Psi:M\lto M$ 
such that 
$$
\int_M c_T(x,\Psi(x))d\mu_0
\leq 
\int_M c_T(x,F(x))d\mu_0
$$
for any transport map $F$.
It turns out that, under the assumption that 
$\mu_0$ has no atoms, a transport map is optimal if and only if
the induced transport plan is an optimal transport plan,
see  \cite{Am:00}, Theorem 2.1.
In other words, we have 
$$
\inf _{\Psi} \int_M c(x, \Psi(x))d\mu_0(x)
=C(\mu_0,\mu_1),
$$
where the 
infimum is taken on the set of transport maps
from $\mu_0$ to $\mu_1$.
This is a general result which holds for any  continuous cost $c$.
It is a standard question, which turns out to be very
hard for certain cost functions,
whether the infimum above is reached, or in other words
whether there exists an optimal transport plan which is induced
from a transport map.
Part of the result below is that this holds true
in the case of the cost $c_0^T$. The method we use
to prove this is an elaboration on ideas due to Brenier, see \cite{Br:91}
and developed for instance in \cite{GaMc:96}, (see also \cite{Ga:ha})
and \cite{Ca:un}, which is certainly the closest to our needs.\vs

\noindent
\textbf{Theorem B.}
\begin{itshape}
Assume that $\mu_0$ is absolutely continuous with respect to the
Lebesgue class on $M$.
Then for each final measure $\mu_T$, there exists one and only one
interpolation $\mu_t,t\in[0,T]$, and each interpolating measure
$\mu_t,t<T$ is absolutely continuous.
In addition, there exists a family $\Psi_0^t: M\rightarrow M$, $t\in ]0,T]$,
of Borel maps such that $(Id\times \Psi_0^t)_{\sharp}\mu_0$
is the only optimal transfer plan in $\mK(\mu_0,\mu_t)$
for the cost function $c_0^t$.
Consequently, we have 
$$
\int_Mc_0^t (x,\Psi_0^t(x))d\mu_0(x)=C_0^t(\mu_0,\mu_t),~~0<t\leq T.
$$
If $\mu_T$, instead of $\mu_0$, is assumed absolutely continuous,
then there exists one and only one interpolation, and each interpolating 
measure $\mu_t,t\in ]0,T]$ is absolutely continuous. 
\end{itshape}\vs

This theorem will be proved and commented in section \ref{AC}.
\subsection{Mather theory}

Let us now assume that the Lagrangian function
is defined for all times, $L\in C^2(TM\times\Rm,\Rm)$
and, in addition to the standing hypotheses, satisfies the
periodicity condition 
$$
L(x,v,t+1)=L(x,v,t)
$$
for all $(x,v,t)\in TM\times\Rm$.
A Mather  measure, see \cite{Ma:91}, is a compactly supported 
probability measure $m_0$ on $ TM $
which is invariant in the sense that $(\psi_0^1)_{\sharp}m_0=m_0$
and is minimizing the action 
$$
A_0^1(m_0)=\int_{TM\times[0,1]} L(\psi_0^t(x,v),t)dm_0 dt.
$$
The major discovery  of  \cite{Ma:91} is  that Mather
measures are supported on the graph of a Lipschitz vectorfield.
Let us call $\alpha$ the action of Mather measures
--this number is the value at zero
of the $\alpha$ function defined by Mather in \cite{Ma:91}.
Let us now explain how this theory of Mather is related to,
and can be recovered from, the content of our paper.\vs

\noindent
\textbf{Theorem C.}
\begin{itshape}
We have 
$$
\alpha
=\min _{\mu} C_0^1(\mu,\mu),
$$
where the minimum is taken on the set of probability measures on $M$.
The mapping 
$m_0\lmto \pi_{\sharp}m_0$ is a bijection
between the set of Mather measures $m_0$ and the set
of probability measures $\mu$ on $M$ satisfying 
$C_0^1(\mu,\mu)=\alpha$.
There exists a Lipschitz vectorfield $X_0$ on $M$
such that all the Mather measures are supported on the graph of $X_0$.
\end{itshape}
\vs

This theorem will be proved in section \ref{AMT},
where the bijection between Mather measures
and measures minimizing $C_0^1(\mu, \mu)$ will
be precised.

\section{Direct variational problems}\label{problems}
We state two different variational
problems whose solutions are the interpolated transports.
We believe that these problems are interesting on their
own. They will also be used to prove Theorem A.

\subsection{Measures}\label{measures}
This formulation parallels Mather's theory.
It can also be related to the generalized curves of L. C. Young.
Let $\mu_0$ and $\mu_T$ be two probability Borel measures on $M$.
Let $m_0\in \mB_1(TM)$ be a    Borel probability
measure on the tangent bundle
$TM$.
We say that $m_0$ is an initial transport measure
if the measure $\eta$ on $M\times M$ given by
$$
\eta=(\pi \times (\pi\circ \psi_0^T))_{\sharp}m_0
$$
is a transport plan, where $\pi:TM\lto M$ is the canonical projection.
We call $\mI(\mu_0,\mu_T)$ the set of initial transport
measures.
To an initial transport measure $m_0$,
we associate the continuous family of measures
$$
m_t=(\psi_0^t)_{\sharp}m_0, t\in[0,T]
$$
on $TM$,
and the measure $m$ on $TM\times [0,T]$
given by
$$
m=m_t \otimes  dt=\big((\psi_0^t)_{\sharp} m_0\big)\otimes  dt.
$$
Note that the linear mapping 
$m_0\lmto m=((\psi_0^t)_{\sharp} m_0) \otimes  dt$
is continuous from $\mB(TM)$ to $\mB(TM\times[0,T])$ endowed with the
weak topology, see appendix.
\begin{lem}
The measure $m$ satisfies the relation
\begin{equation}\label{closed}
\int_{TM\times[0,T]} \partial_tf(x,t)+\partial_xf(x,t)\cdot v\,  dm(x,v,t)
=\int _M f_T d\mu_T-\int _Mf_0 d\mu_0
\end{equation}
for each function $f\in C^1(M\times[0,T],\Rm)$,
where $f_t$ denotes the function $x\lmto f(x,t)$.
\end{lem}

\proof
Setting $\widetilde f(x,v,t)=f(x,t)$,
$g_1(x,v,t)=\partial_t f(x,t)=\partial_t \widetilde f(x,v,t)$ and
$g_2(x,v,t)=\partial_x f(x,t)\cdot v$,
we have 
$$\int_{TM\times[0,T]} \partial_tf(x,t)+\partial_xf(x,t)\cdot v\,   dm(x,v,t)
=\int_0^T \int _{TM} 
(g_1 + g_2)\circ \psi_0^t\,  dm_0 dt. 
$$
Noticing that, in view of equation (\ref{elflow}), we have
$$
\partial_t (\widetilde f \circ \psi_0^t)
= g_1 \circ  \psi_0^t+g_2 \circ \psi_0^t
$$
we obtain
that 
$$
\int_{TM\times[0,T]} \partial_tf(x,t)+\partial_xf(x,t)\cdot v \, dm(x,v,t)
=\int_{TM} (\widetilde f \circ\psi_0^T-\widetilde f)dm_0
=\int _M f_T d\mu_T-\int _Mf_0 d\mu_0
$$
as desired. 
\qed%
\begin{defn}\label{defn invariant}
A finite Borel measure on $TM\times[0,T]$ which satisfies 
(\ref{closed}) is called a transport measure.
We denote by $\mM(\mu_0,\mu_T)$ the set of transport measures.
A transport measure which is induced 
from an initial measure $m_0$ is called an invariant transport 
measure.
The action of the transport measure $m$
is defined by 
$$
A(m)=\int_{TM\times[0,T]}
L(x,v,t)dm\in \Rm\cup \{\infty\}
$$
The action $A(m_0)$ of an initial transport measure
is defined as the action of the associated transport measure $m$.
We will also denote this action by $A_0^T(m_0)$ when
we want to insist on the time interval.
We have
$$
A_0^T(m_0)=\int_{TM\times[0,T]} L(\psi_0^t(x,v),t)dm_0 dt.
$$
\end{defn}
 
Notice that initial tranport measures exist:
\begin{prop}\label{planmesure}
The mapping 
$\big(\pi\times (\pi\circ \psi_0^T)\big)_{\sharp}:\mI(\mu_0,\mu_T)
\lto \mK(\mu_0,\mu_T)$
is surjective. In addition,
for each transport plan $\eta$, there exists a
compactly supported  initial
transport measure $m_0$ such that
$(\pi\times (\pi\circ \psi_0^T))_{\sharp}m_0=\eta$
and such that 
$$
A(m_0)=\int _{M\times M}c_0^T(x,y)d\eta. 
$$
\end{prop}
\proof
By Proposition \ref{extremals},
there exists a compact set $K\in TM$ such that
if $\gamma(t):[0,T]\lto M$ is a minimizing extremal,
then the lifting $(\gamma(t),\dot \gamma(t))$
is contained in $K$ for each $t\in[0,T]$.
We shall prove that, for each probability measure $\eta \in \mB(M\times M)$,
there exists a probability  measure $m_0\in \mB(K)$ such that
$(\pi\times (\pi\circ \psi_0^T))_{\sharp}m_0=\eta$
and such that 
$$
A(m_0)=\int _{M\times M}c_0^T(x,y)d\eta. 
$$
Observing that 
\begin{itemize}
\item the mappings $m_0\lto (\pi\times (\pi\circ \psi_0^T))_{\sharp}m_0$ 
and
$m_0\lmto A(m_0)$ are linear and continuous on 
the space $\mB_1(K)$ of probability measures supported on $K$,
\item the set $\mB_1(K)$ is compact for the weak topology,
and the action $A$ is continuous on this set,
\item the set of probability measures on $M\times M$ 
is the compact convex  closure of the set of Dirac  probability measures
(probability measures supported in one point), see e. g. 
\cite{Bi:99}, p. 73,
\end{itemize}
it is enough to prove the result when $\eta$ is a Dirac probability measure
(or equivalently 
when $\mu_0$ and $\mu_T$ are Dirac probability measures).
Let  $\eta$ be  the Dirac probability measure supported
at $(x_0,x_1)\in M\times M$.
Let $\gamma(t):[0,T]\lto M$ be a minimizing extremal
with boundary conditions $\gamma(0)=x_0$
and $\gamma(T)=x_1$.
In view of the choice of $K$, we have $(\gamma(0),\dot \gamma(0))\in K$. 
Let $m_0$ be the Dirac probability measure supported
at $(\gamma(0),\dot\gamma(0))$.
It is straightforward that 
$m_t$ is then the Dirac measure supported at $(\gamma(t),\dot \gamma(t))$,
so that 
$$
A(m_0)=\int_0^T L dm_t dt
=\int _0^T L(\gamma(t),\dot \gamma(t),t) dt
=c_0^T(x_0,x_1)
=\int_{M\times M} c_0^T d\eta
$$
and 
$$
(\pi \times (\pi \circ \psi_0^T))_{\sharp}
m_0=\eta.
$$
\qed

Although we are going to build minimizers by other means,
we believe the following result is worth being mentioned.

\begin{lem}
For each real number $a$, the set 
$\mM^a(\mu_0,\mu_T)$ of transport measures $m$  which satisfy
$A(m)\leq a$, as well as the set
$\mI^a(\mu_0,\mu_T)$ of initial  transport  measures $m_0$  which satisfy
$A_0^T(m_0)\leq a$,
are compact.
As a consequence, there exist optimal initial transport measures,
and optimal transport measures.
\end{lem}
\proof
This is an easy application of the  Prohorov theorem, see 
the Appendix.
\qed

Now we have seen that the problem of finding optimal transport
measures is well-posed, let us describe its solutions.
\begin{thm}
We have 
$$
C_0^T(\mu_0,\mu_T)=\min _{m\in \mM(\mu_0,\mu_T)}
A(m)=\min _{m_0\in \mI(\mu_0,\mu_T)}A(m_0).
$$
The mapping 
$$
m_0\lmto m=\big((\psi_0^t)_{\sharp}m_0\big)\otimes dt
$$
between the set $\mO\mI$ of optimal initial measures
and the set $\mO\mM$ of optimal transport measures
is a bijection.
There exists a bounded and locally Lipschitz vectorfield
$X(x,t):M\times]0,T[\lto TM$
such that, for each optimal initial measure 
$m_0\in \mO\mI$, the measure 
$m_t=(\psi_0^t)_{\sharp}m_0$
is supported on the graph of $X_t$ for each $t\in ]0,T[$.
\end{thm}

The proof will be given in section \ref{COM}.
Let us just notice now that the inequalities 
$$
C_0^T(\mu_0,\mu_T)
\geq  \min _{m_0\in \mI(\mu_0,\mu_T)}A(m_0)
\geq \min _{m\in \mM(\mu_0,\mu_T)}
A(m)
$$
hold in view of Proposition \ref{planmesure}.

\subsection{Currents}\label{currents}
This formulation finds its roots on one hand 
in the works  of Benamou and Brenier, see \cite{BeBr}, and then Brenier,
see \cite{Br:00},
and on the other hand in the work of Bangert \cite{Ba:99}.
Let $\Omega^0(M\times[0,T])$ be the set of continuous one-forms 
on $M\times [0,T]$, endowed with the uniform norm.
We will often decompose forms $\omega \in \Omega^0(M\times[0,T])$
as 
$$\omega =\omega^x + \omega^t dt,
$$
where
$\omega^x$ is a time-dependent form on $M$ and $\omega^t$
is a continuous function on $M\times[0,T]$.
To each continuous linear form $\chi$ on $\Omega^0(M\times [0,T])$,
we associate its time component 
$\mu_{\chi}$, which is the measure on $M\times[0,T]$ defined by
$$
\int _{M\times [0,T]}fd\mu_{\chi}=\chi(fdt)
$$
for each continuous function $f$ on $M\times [0,T]$.
A Transport current between $\mu_0$ and $\mu_T$
is a continuous linear  form 
$\chi$
on $\Omega^0(M\times[0,T])$
which satisfies the two conditions:
\begin{enumerate}
\item 
The measure $\mu_{\chi}$ is  non negative (and bounded). 
\item 
$d\chi=\mu_T\otimes \delta_T  -
\mu_0\otimes \delta_0$,
which means that 
$$
\chi(df)=\int_M f_Td\mu_T-\int_Mf_0d\mu_0
$$
for each smooth (or equivalently $C^1$) function
$f:M\times [0,T]\lto \Rm$.
\end{enumerate}
We call $\mC(\mu_0,\mu_T)$
the set of transport currents from $\mu_0$
to $\mu_T$.
The set  $\mC(\mu_0,\mu_T)$ is a closed  convex subset
of   $\big[\Omega^0(M\times[0,T])\big]^*$.
We will endow  $\mC(\mu_0,\mu_T)$ with the weak topology
obtained as the restriction of the weak-$*$ topology of 
 $\big[\Omega^0(M\times[0,T])\big]^*$.
Transport currents should be thought of as vectorfields
whose components are measures, the last component being $\mu_{\chi}$. 

If $Z$ is a  bounded measurable 
vectorfield on $ M \times [0,T]$,
and if $\nu$ is a finite non-negative measure on $M\times[0,T]$,
we define the current $Z\wedge \nu$
by 
$$
Z\wedge \nu(\omega):=\int_{M\times[0,T]}
\omega(Z)d\nu.
$$
Every transport current can be written  in this way, see 
\cite{Fe:69} or \cite{GiMoSo}.
As a consequence, currents extend as linear forms
on the set $\Omega_{\infty}(M\times[0,T])$
of bounded measurable one-forms.
If $I$ is a Borel subset of the interval $[0,T]$,
it is therefore possible to define the restriction 
$\chi_I$ of the current $\chi$ to $I$
by the formula
$
\chi_I(\omega)=\chi(1_I\omega)
$, where $1_I$ is the indicatrix of $I$.

\begin{lem}\label{desint}
If $\chi$ is a transport current,
then 
$$
\tau_{\sharp} \mu_{\chi}
=dt,
$$
where $\tau$ is the projection onto $[0,T]$, see appendix.
As a consequence, there exists a measurable familty $\mu_t,t\in]0,T[$
of probability measures on $M$ such that
$\mu_{\chi}=\mu_t\otimes dt$, see appendix.
There exists a set $I\subset]0,T[$ of total measure
such that the relation 
\begin{equation}\label{difference}
\int_Mf_td\mu_t=\int_Mf_0d\mu_0 +\chi_{[0,t[}(df)
\end{equation}
holds 
for each $C^1$ function $f:M\times [0,T]\lto M$
and each $t\in I$.
\end{lem}
\proof
Let $g:[0,T]\lto \Rm$ be a continuous function.
Setting  $G(t)=\int_0^t g(s)ds$,  we  observe that
$$\int_{M\times[0,T]}g d\mu_{\chi}
=\chi(dG)=\int_MG_Td\mu_T-\int_M G_0d\mu_0
=G(T)-G(0)=\int_0^Tg(s)ds.
$$
This implies that 
$
\tau_{\sharp} \mu_{\chi}
=dt.
$
As a consequence, the measure $\mu_{\chi}$
can be desintegrated as $\mu_{\chi}= \mu_t\otimes dt$.
We claim that, for each $C^1$ function $f:M\times [0,T]\lto M$,
the relation (\ref{difference})
holds for almost every $t$.
Since the space $C^1(M\times [0,T],\Rm)$ is
separable, the claim implies the existence of a set $I\subset ]0,T[$
of full Lebesgue measure such that (\ref{difference})  holds
for all $t\in I$ and all $f\in C^1(M\times [0,T],\Rm)$.
In order to prove the claim, 
let us fix a function $f$ in  $C^1(M\times [0,T],\Rm)$.
For each function 
$g\in C^1([0,T],\Rm)$,
we have
$$
\chi(d(gf))=\chi(g'fdt)+\chi(gdf)
$$
hence
$$
g(T)\int_M f_Td\mu_T-g(0)\int_Mf_0d\mu_0
=
\int _0^T g'(t)\int_Mf_td\mu_t dt
+\chi(gdf).
$$
By applying this relation to a sequence of $C^1$
functions $g$ approximating $1_{[0,t[}$, we get, at the limit
$$
-\int_Mf_0d\mu_0
=
-\int_Mf_td\mu_t
+\chi_{[0,t[}(df)
$$
at every Lebesgue point of the function $t\rightarrow \int_M f_td\mu_t$ .
\qed

If $\mu_0=\mu_T$, an easy example of transport current is given by 
$\chi(\omega)=\int_{M}\int_{0}^T \omega^t dtd\mu_0$.
Here are some more interesting  examples.

\noindent
\textbf{Regular transport currents.}
The transport current $\chi$ is called regular
if there exists a 
bounded measurable section $X$ of the projection 
$ TM\times[0,T]\lto M\times [0,T]$,
and a  a non-negative measure $\mu$ on $M\times[0,T]$
such that $\chi=(X,1)\wedge \mu$.
The time component
of the current $(X,1)\wedge \mu $ is $\mu$.
In addition, if $(X,1)\wedge \mu=(X',1)\wedge \mu$
for two  vectorfields $X$ and $X'$, then $X$ and $X'$
agree $\mu$-almost everywhere.

\begin{itshape}
The current $\chi=(X,1)\wedge \mu$, with $X$ bounded,
is a regular transport current if and only if 
there exists a (unique)  continuous family 
$\mu_t\in \mB_1(M), t\in [0,T]$ 
(where $\mu_0$ and $\mu_T$ are the 
transported measures)
such that $\mu_{\chi}=\mu_t\otimes dt$
and such that the transport equation 
$$
\partial_t \mu_t +\partial_x.(X\mu_t)=0
$$
holds in the sense of distributions on $M\times ]0,T[$.
The relation
$$
\int_M f_td\mu_t-\int_M f_sd\mu_s
= \chi_{[s,t[}(df)
$$
then holds for each $C^1$ function $f$ and each $s\leq t$
in $[0,T]$.
\end{itshape}

In order to prove that the family 
$\mu_t$ can be chosen  continuous,
pick a function $f\in C^1(M,\Rm)$ and
 notice that the equation
$$
\int_Mfd\mu_t-\int_Mfd\mu_s
=\chi_{[s,t[}
(df)
=
\int_s^t \int_M df\cdot X_{\sigma}d\mu_\sigma d\sigma 
$$
holds for all $s\leq t$ in a subset of total measure $I\subset[0,T]$.
Note that this relation also holds if $s=0$ and $t\in I$ and
if $s\in I$ and  $t=T$.
Since the function $\sigma \lmto \int_M df\cdot X_{\sigma}d\mu_{\sigma}$
is bounded, we conclude that the function
$t\lmto \int_Mfd\mu_t$ is Lipschitz on $I\cup \{0,T\}$ for each 
$f\in C^1(M,\Rm)$, with a Lipschitz constant which
depends only on
$\|df\|_\infty\cdot\|X\|_{\infty}$.
The family $\mu_t$ is then Lipschitz on $I\cup \{0,T\}$ for 
the  $1$-Wasserstein 
distance on probability measures, see \cite{Vi:03,Du:02,AmGiSa} 
for example, the Lipschitz constant depending
only on $\|X\|_{\infty}$. It suffices to remember that, on the compact
manifold $M$, the $1$-Wasserstein distance on probabilities
is topologically equivalent to the weak topology, see for example
\cite{Yo}, (48.5)
or \cite{Vi:03}.

\noindent
\textbf{Smooth transport currents.}
A regular transport current is said smooth if it can be written on the form
$(X,1)\wedge \lambda$ with a bounded 
vectorfield $X$ smooth on $M\times]0,T[$  and
a measure $\lambda$
that   has a positive smooth
density with respect to the Lebesgue class in any chart in $M\times]0,T[$.
Every   transport  current in $\mC(\mu_0,\mu_T)$ 
can be approximated by smooth
transport currents, but we shall not use such
approximations.

\noindent
\textbf{Lipschitz regular transport currents.}
A regular
transport current is said Lipschitz regular
if it can be written
in the form $(X,1)\wedge \mu $
with  a vectorfield $X$
which is bounded and locally  Lipschitz on $M\times]0,T[$.
Smooth currents are Lipschitz regular.
Lipschitz regular transport currents 
have a remarkable structure:

\begin{itshape}
If $\chi=(X,1)\wedge \mu $ is a Lipschitz regular transport current
with $X$ bounded and  locally Lipschitz on $M\times]0,T[$, then
 we have
$$(\Psi_s^t)_{\sharp}\mu_s=\mu_t$$
where $\Psi_s^t, (s,t)\in ]0,T[^2$, denotes the flow of
the Lipschitz vectorfield  $X$
from time $s$ to time $t$,
and $\mu_t$ is the continuous family of
probability measures such that $\mu_{\chi}=\mu_t\otimes dt$.
\end{itshape}

This statement follows from standard representation
results for solutions of the transport equation,
see for example \cite{Am:transport} or \cite{AmGiSa}.

\noindent
\textbf{Transport current induced from a transport measure.}
To a transport measure $m$, we associate the transport current
$\chi_m$ defined by
$$
\chi_m(\omega)=\int_{TM\times[0,T]}\big(\omega^x(x,t)\cdot v
+\omega^t(x,t)\big)dm(x,v,t)
$$
where the form $\omega$
is decomposed as $\omega =\omega^x+\omega^t dt$.
Note  that the time component of the current $\chi_m$ is $\pi_{\sharp}m$.
We will see in Lemma \ref{inequalities} that 
$$
A(\chi_m)\leq A(m)
$$
with the following definition of the action $A(\chi)$
of a current, with equality if $m$ is concentrated on the
graph of any bounded  vectorfield $M\times[0,T]\rightarrow TM$.

\begin{lem}
For each transport current $\chi$,
the numbers
\begin{align*}
A_1(\chi) &=
\sup _{\omega\in \Omega^0} \Big(\chi(\omega^x,0)
-\int _{M\times[0,T]}
H(x,\omega^x(x,t),t)
d\mu_{\chi}\Big)\\
A_2(\chi) &=\sup _{\omega\in \Omega^0} \Big(\chi(\omega)
-\int _{M\times[0,T]}
\big(H(x,\omega^x(x,t),t)+\omega^t\big)
d\mu_{\chi}\Big)\\
A_3(\chi) &=\sup _{\omega\in \Omega^0} \Big(\chi(\omega)
-T\sup _{(x,t)\in M\times[0,T]}
\big(H(x,\omega^x(x,t),t)+\omega^t\big)\Big)\\
A_4(\chi) &=\sup _{\omega\in \Omega^0;\omega^t+H(x,\omega^x,t)\leq 0}
\chi(\omega)\\
A_5(\chi) &=\sup _{\omega\in \Omega^0;\omega^t+H(x,\omega^x,t)\equiv 0}
\chi(\omega),
\end{align*}
are  equal.
In addition the number $A_i^{\infty}(\chi)$ obtained by replacing
in the above suprema the set $\Omega^0$ of continuous forms
by the set $\Omega_{\infty}$ of bounded measurable forms 
also have the same value.
\end{lem}
The last remark in the statement has been added in the last
version of the paper and is inspired from \cite{PaStGr}.

\proof
It is straightforward
that $A_1=A_2$, this just amounts to
simplifying the term $\int \omega^td\mu_{\chi}$.
Since $\mu_{\chi}$
is a non-negative measure which satisfies
$\int_{M\times[0,T]} 1 d\mu_{\chi}=T$, we have
$$
\int _{M\times[0,T]}
\big(H(x,\omega^x(x,t),t)+\omega^t\big)d\mu_{\chi}
\leq T\sup _{(x,t)\in M\times[0,T]}
\big(H(x,\omega^x(x,t),t)+\omega^t\big)
$$
so that $A_3(\chi)\leq A_2(\chi)$.
In addition, we obviously have $A_5(\chi)\leq A_4(\chi) \leq A_3(\chi)$.
Now notice, in $A_2$, that the quantity 
$$
\chi(\omega)-\int _{M\times[0,T]}
\big(H(x,\omega^x(x,t),t)+\omega^t\big)
d\mu_{\chi}
$$
does not depend on $\omega^t$.
Let us consider the form
$\tilde \omega=(\omega^x,-H(x,\omega^x,t))$,
which satisfies the equality $H(x,\tilde \omega^x,t)+\tilde \omega^t\equiv 0$.
We get, for each form $\omega$,
$$
\chi(\omega^x,0)-\int _{M\times[0,T]}
H(x,\omega^x(x,t),t)
d\mu_{\chi}=\chi(\tilde \omega)\leq A_5(\chi)
$$
Hence $A_1(\chi)\leq A_5(\chi)$.
Exactly the same proof shows that the numbers $A_i^{\infty}(\chi)$
are equal.
In order to end the proof, it is enough to check that 
$A_2(\chi)=A_2^{\infty}(\chi)$.
Writing the current $\chi$ on the form $Z\wedge \nu$ with a bounded
vectorfield $Z$ and a measure $\nu\in \mB_+(M\times[0,T])$,
we have
$$
A_2(\chi)
=
\sup_{\omega\in \Omega^0}\Big(
\int_{M\times [0,T]}
\omega(Z)d\nu-\int _{M\times[0,T]}
\big(H(x,\omega^x(x,t),t)+\omega^t\big)
d\mu_{\chi}\Big)
$$
and 
$$
A_2^{\infty}(\chi)=
\sup_{\omega\in \Omega_{\infty}}\Big(
\int_{M\times [0,T]}
\omega(Z)d\nu-\int _{M\times[0,T]}
\big(H(x,\omega^x(x,t),t)+\omega^t\big)
d\mu_{\chi}\Big).
$$
The desired result follows by density 
of continuous functions  in $L^1(\nu+\mu_{\chi})$.
\qed

\begin{defn}
We denote by $A(\chi)$
and call action of the transport current $\chi$
the common value of the numbers $A_i(\chi)$.
\end{defn}

The existence of currents of finite action 
follows from the following:

\begin{lem}\label{inequalities}
We have 
$$
A(\chi)
=\int_{M\times[0,T]}
L(x,X(x,t),t)d\mu
$$
for each  regular current $\chi = (X,1)\wedge \mu$.
If $m$ is a transport measure, and if $\chi_m$ is the associated
transport current, then $A(\chi_m)\leq A(m)$, with equality if $m$
is supported on the graph of a bounded
Borel vectorfield.
As a consequence, we have the inequalities 
$$
C_0^T(\mu_0,\mu_T)\geq 
\min_{m_0\in \mI(\mu_0,\mu_T)} A(m_0)
\geq
\min_{m\in \mM(\mu_0,\mu_T)} A(m)
\geq 
\min_{\chi\in \mC(\mu_0,\mu_T)} A(\chi).
$$
\end{lem}
\proof
For each bounded measurable  form $\omega$,
we have 
$$
\int _{M\times[0,T]}
\omega^x(X)
-H(x,\omega^x(x,t),t)
d\mu
\leq 
\int _{M\times[0,T]}
L(x,X(x,t),t)
d\mu,
$$
so that 
$$
A((X,1)\wedge \mu)\leq \int _{M\times[0,T]}
L(x,X(x,t),t)
d\mu.
$$
On the other hand, taking the form 
$\omega_0^x(x,t)=\partial_v L(x,X(x,t),t)$ 
we obtain the pointwise equality
$$
L(x,X(x,t),t)=\omega^x_0(X)-H(x,\omega^x_0(x,t),t)
$$
and by integration
$$
\int _{M\times[0,T]}
L(x,X(x,t),t)
d\mu 
=
\int _{M\times[0,T]}
\omega_0^x(X)
-H(x,\omega_0^x(x,t),t)
d\mu
\leq A((X,1)\wedge \mu).
$$
This ends the proof of the equality of the two forms
of the action of regular  currents. 
Now if $\chi_m$ is the current associated to a transport measure $m$,
then we have, for each bounded form 
$\omega\in \Omega^0(M\times[0,T])$,
$$
\chi_m(\omega)-
\int _{M\times[0,T]}
\omega^t(x,t)+H(x,\omega^x(x,t),t)
d\mu_{\chi}
=
\int _{TM\times[0,T]}
\omega^x(v)-H(x,\omega^x(x,t),t)
dm
$$
by definition of $\chi_m$, so that 
$$
A(\chi_m)\leq\int _{TM\times[0,T]}
L(x,v,t)
dm=
 A(m)
$$
by the Legendre inequality.
In addition, if there exists a bounded measurable
vectorfield $X:M\times[0,T]\lto TM$ such
that the graph of $X\times\tau$ supports $m$, then we can consider
the form  $\omega_0^x$ associated to $X$ as above, 
and we get the equality for this form.
\qed

Although we are going to provide explicitly a minimum of
$A$, we believe the following Lemma is worth being mentioned.

\begin{lem}
The functional $A:\mC(\mu_0,\mu_T)\lto \Rm\cup\{+\infty\}$
is convex and lower semi-continuous, both for the strong and weak-$^*$
topologies on $[\Omega^0(M\times[0,T])]^*$. Moreover it is coercive
with respect to the strong topology and
hence it has a minimum.
\end{lem}
\proof
First note that $A(\chi)<\infty$ if $\chi$ is the transport current
corresponding to an initial transport measure in $\mM(\mu_0,\mu_T)$
arising from a transport plan.
Let us define the continuous 
convex function $\Hm_T:\Omega^0(M\times[0,T])\lto \Rm$
by
$$
\Hm_T(\omega)=T\sup_{(x,t)\in M\times[0,T]}
H(x,\omega^x(x,t),t)+\omega^t.
$$
Then the action is the restriction to $\mC(\mu_0,\mu_T)$
of the Fenchel conjugate 
$A=\Hm^*:[\Omega^0(M\times[0,T])]^*\lto \Rm\cup\{+\infty\}$.
In other words,
 $A$ is the supremum over $\omega$ of the family of affine functionals
$$\chi\rightarrow \chi(\omega)-\Hm_T(\omega)$$
that are continuous both for the strong and weak-$^*$ topologies. Hence
$A$ is convex and lower semi-continuous for both topologies.
Since
$$A(\chi)\geq \sup_{\|\omega\|\leq 1}\chi(\omega)
- \sup_{\|\omega\|\leq 1}\Hm_T(\omega),$$
$A$ is coercive. The existence of a minimizer is standard:
any minimizing sequence $(\chi_n)$ is bounded (thanks to coercivity)
and has a weakly-$*$ convergent subsequence (because $\Omega^0(M\times[0,T])$
is a separable Banach space). By lower semicontinuity, 
its weak-$^*$ limit is a minimizer.
Note that $\mC(\mu_0,\mu_T)$ is weakly$^*$ closed.
\qed

\begin{thm}
We have 
$$
C_0^T(\mu_0,\mu_T)=\min_{\chi\in \mC(\mu_0,\mu_T)} A(\chi)
$$
where the minimum is taken on all transport currents 
from $\mu_0$ to $\mu_T$.
Every optimal transport current is Lipschitz regular.
Let  $\chi=(X,1)\wedge \mu$ be  an optimal
transport current,
with $X$ locally Lipschitz on $M\times]0,T[$.
The measure  $m= (X\times \tau)_{\sharp} \mu\in \mB_+(TM\times]0,T[)$ 
is an optimal transport
measure, and $\chi$ is the transport current induced from $m$. 
Here $\tau:TM\times[0,T]\lto [0,T]$ is the projection on the second 
factor, see appendix.
We have 
$$
C_0^T(\mu_0,\mu_T)=A(m)=A(\chi)=\int_{M\times[0,T]} L(x,X(x,t),t)d\mu_{\chi}.
$$
\end{thm}

This result will be proved in \ref{OTC}
after some essential results on the dual approach have
been established.
\section{Hamilton-Jacobi equation}\label{sectionHJ}
Most of the results stated so far can be proved
by direct approaches using Mather's shortening Lemma,
which in a sense is an improvement
of the initial observation of Monge, see \cite{Ma:91}
and \cite{Ba:99}.
We shall however base our proofs on the use of 
the Hamilton-Jacobi equation,
in the spirit of Fathi's approach to Mather theory, see \cite{Fa:un},
which should be associated to Kantorovich
dual approach of the transportation problem.

\subsection{Viscosity solutions and semi-concave functions}
It is certainly useful to recall the main properties of
viscosity solutions in connection with semi-concave functions.
We  will not  give proofs, and instead refer 
to \cite{Fa:un}, \cite{FaSi},  \cite{CaSi}, as well 
as the appendix in \cite{Be:05}. 
We will consider the Hamilton-Jacobi equation
\begin{equation}\tag{$HJ$}
\partial_tu+H(x,\partial_x u,t)=0.
\end{equation}
The function $u:M\times [0,T]\lto M$
is called $K$-semi-concave if, for each chart $\theta\in \Theta$
(see appendix), the function
$$
(x,t) \lmto u(\theta(x),t)-K(\|x\|^2+t^2)
$$
is concave on $B_3\times [0,T]$.
 The function $u$ is called semi-concave if it is
$K$-semi-concave for some $K$. A function $u:M\times]0,T[\lto M$
is called locally semi-concave if it is semi-concave on each
$M\times[s,t]$, for $0<s<t< T$.
The following regularity result follows from Fathi's work,
see \cite{Fa:un} and also \cite{Be:05}.
  
\begin{prop}\label{Fathilip}
Let $u_1$ and $u_2$ be two $K$-semi-concave functions.
Let $A$ be the set of minima of  the function $u_1+u_2$.
Then the functions $u_1$ and $u_2$ are differentiable on $A$,
and $du_1(x,t)+du_2(x,t)=0$ at each point of $(x,t)\in A$.
In addition, the mapping $du_1:M\times [0,T]\lto T^*M$
is $CK$-Lipschitz continuous on $A$, where $C$ is a universal constant.
\end{prop}

\begin{defn}
We say that the function $u:M\times]s,t[\lto \Rm$
is a viscosity solution of $(HJ)$ if the equality
$$
u(x,\sigma)=\min_{y\in M} u(y,\zeta)+c_{\zeta}^{\sigma}(y,x)
$$ 
holds for all $x\in M$ and all $s<\zeta<\sigma<t$.

We say that the function $\breve u:M\times]s,t[\lto \Rm$
is a backward viscosity solution of $(HJ)$ if the equality
$$
\breve u(x,\sigma)=\max_{y\in M} \breve u(y,\zeta)-c^{\zeta}_{\sigma}(x,y)
$$ 
holds for all $x\in M$ and all $s<\sigma <\zeta<t$.

We say that the function $v:M\times]s,t[\lto \Rm$
is a viscosity sub-solution of (HJ)
if the inequality
$$
v(x,\sigma)\leq  v(y,\zeta)+c_{\zeta}^{\sigma}(y,x)
$$
holds for all $x$ and $y$ in $ M$ and all $s<\zeta<\sigma<t$.

Finally we will  say that the function $v:M\times[s,t]\lto \Rm$ is a 
continuous viscosity solution 
(subsolution, backward solution)
of (HJ) 
if it is continuous on $M\times[s,t]$ and if
if $v|_{M\times]s,t[}$ is a viscosity solution of (HJ) 
(subsolution, backward solution).
\end{defn}

Notice that both viscosity solutions
and backward viscosity solutions are viscosity sub-solutions.
That these definitions are equivalent 
in our setting
to the usual ones
is studied in the references listed above, but is not 
useful for our discussion.
The only fact which will be used is that, for a $C^1$
function $u:M\times ]s,t[\lto \Rm$, being a viscosity
solution (or a backward viscosity solution)
is equivalent to being a pointwise solution of $(HJ)$,
and being a viscosity sub-solution is equivalent to satisfy 
pointwise the inequality
$\partial_tu+H(x,\partial_xu,t)\leq 0$.

\noindent
\textbf{Differentiability  of viscosity solutions.}
Let $u\in C(M\times[0,T[,\Rm)$ be a viscosity solution of $(HJ)$
(on the interval $]0,T[$).
We have the expression 
$$
u(x,t)=\min_{\gamma}
u(\gamma(0),0)+\int_0^tL(\gamma(\sigma),\dot \gamma(\sigma),\sigma)d\sigma
$$
where the minimum is taken on the set 
of curves $\gamma\in C^2([s,t],M)$
which satisfy the final condition $\gamma(t)=x$.
Let us denote by $\Gamma(x,t)$ the set of minimizing curves
in this expression, which are obviously minimizing
extremals of $L$. 
We say that $p\in T^*_xM$ is a proximal super-differential
of the function $u:M\lto \Rm$ at point $x$ if there 
exists a smooth function $f:M\lto \Rm$
such that $f-u$ has a minimum at $x$ and $d_xf=p$.

\begin{prop}\label{proximal}
Let us fix a point $(x,t)\in M\times]0,T[$.
The function $u_t$ is differentiable at $x$
if and only if 
the set
$\Gamma(x,t)$ contains a single element
$\gamma$, and then $\partial_xu(x,t)=\partial_vL(x,\dot \gamma(t),t)$.

For all $(x,t)\in M\times]0,T[$ and $\gamma \in \Gamma(x,t)$,
we set $p(s)=\partial_vL(\gamma(s), \dot \gamma(s),s)$.
Then $p(0)$ is a proximal sub-differential of $u_0$
at $\gamma(0)$, and $p(t)$ is a proximal super-differential
of $u_t$ at $x$.
\end{prop}

We finish with an important property on regularity of viscosity solutions:

\begin{prop}
For each continuous function $u_0:M\lto \Rm$,
the viscosity solution 
$$
u(x,t):=\min_{y\in M} u_0(y)+c_0^t(y,x)
$$
is locally semi-concave on $]0,T]$.
If in addition the initial condition $u_0$ is 
Lipschitz, then  $u$ is Lipschitz  on $[0,T]$.

For each continuous function $u_T:M\lto \Rm$,
the viscosity solution 
$$
\breve u(x,t):=\max_{y\in M} u_T(y)-c_t^T(x,y)
$$
is locally semi-convex on $[0,T[$.
If in addition the final condition $u_T$ is 
Lipschitz, then  $u$ is Lipschitz on $[0,T]$.
\end{prop}
\proof
The part concerning semi-concavity of $u$ is proved in \cite{CaSi},
for example.
It implies that the function $u$ is locally Lipschitz on 
 $]0,T]$, hence differentiable almost everywhere.
In addition, at each point of differentiability of 
$u$,
we have $\partial_t u +H(x,\partial_xu,t)=0$
and
 $\partial_xu(x,t)=p(t)=\partial_vL(x,\dot \gamma(t),t)$,
where $\gamma:[0,t]\lto M$ is the only curve 
in $\Gamma(x,t)$.
In order to prove that the function $u$ is Lipschitz,
it is enough to prove that there exists a uniform bound 
on $|p(t)|$.
It is known, see Proposition \ref{proximal},
that $p(0):=\partial_vL(\gamma(0),\dot \gamma(0),0)$
is  a proximal sub-differential of the function $u_0$
at point $\gamma(0)$.
If $u_0$ is Lipschitz, its sub-differentials are bounded:
There exists a constant $K$ such that 
$|p(0)|\leq K$. By completeness, there exists a constant
$K'$,
which depends only on the Lipschitz constant of $u_0$, 
such that $|p(s)|\leq K'$ for all $s\in [0,t]$.
This proves that the function $u$ is Lipschitz.
The statements concerning $\breve u$ are proved in a similar way.
\qed

\subsection{Viscosity solutions and Kantorovich optimal pairs}
Given a Kantorovich optimal pair $(\phi_0,\phi_1)$,
we define the viscosity solution 
$$
u(x,t):= \min_{y\in M} \phi_0(x)+c_0^t(y,x)
$$
and the backward viscosity solution
$$
\breve u(x,t):= \max_{y\in M} \phi_1(y)-c_t^T(x,y)
$$
which satisfy 
$u_0=\breve u_0=\phi_0$,
and $u_T=\breve u_T=\phi_1$.
Note that both $\phi_1$ and $-\phi_0$ are semi-concave
hence Lipschitz, that the function $u$ is
Lipschitz and locally semi-concave on $]0,T]$,
and that the function $\breve u$ is Lipschitz 
and locally semi-convex on $[0,T[$.

\begin{prop}
We have 
\begin{equation}\label{visc}
C_0^T(\mu_0,\mu_T)
=\max _{u}
\Big(
\int_M u_Td\mu_T-\int_M u_0d\mu_0
\Big),
\end{equation}
where the minimum is taken on the set
of continuous viscosity solutions $u:M\times [0,T]\lto \Rm$
of the Hamilton-Jacobi equation $(HJ)$.
The same conclusion holds if the maximum is taken on the set
of continuous backward viscosity solutions.
The same conclusion also holds if the maximum is taken on the set
of continuous viscosity sub-solutions of $(HJ)$.
\end{prop}

\proof
If $u(x,t)$ is a continuous viscosity sub-solution of $(HJ)$,
then it satisfies
$$
u_T(x)-u_0(y)\leq c_0^T(y,x)
$$
for each $x$ and $y\in M$, and so, by Kantorovich duality, 
$$
\Big(
\int_M u_Td\mu_T-\int_M u_0d\mu_0
\Big)\leq C_0^T(\mu_0,\mu_T).
$$
The converse inequality is obtained by using the 
functions $u$ and $\breve u$.
\qed
\begin{defn}\label{dfn: T}
If $(\phi_0,\phi_1)$ is a 
Kantorovich optimal pair,
 then we denote by 
$\mF(\phi_0,\phi_1)\subset C^2([0,T],M)$ the set of curves $\gamma(t)$
such that 
$$
\phi_1(\gamma(T))=\phi_0(\gamma(0))+\int_0^TL(\gamma(t),\dot \gamma(t),t)dt.
$$
We denote by $\mT(\phi_0,\phi_1)\subset M\times]0,T[$
the set
$$
\mT(\phi_0,\phi_1)=\{(\gamma(t),t), t\in ]0,T[, \gamma\in \mF(\phi_0,\phi_1)\}
$$
and by $\tilde \mT(\phi_0,\phi_1)\subset TM\times]0,T[$
the set
$$
\tilde \mT(\phi_0,\phi_1)=
\{(\gamma(t), \dot \gamma(t),t), t\in ]0,T[, \gamma\in \mF(\phi_0,\phi_1)\},
$$ 
which is obviously invariant under the Euler-Lagrange flow.
\end{defn}

\begin{prop}\label{mT}
Let $(\phi_0,\phi_1)$ be a Kantorovich  optimal pair,
and let $u$ and $\breve u$ be the associated viscosity
and backward viscosity solutions.
\begin{enumerate}
\item
We have $\breve u \leq u$, and
$$
\mT(\phi_0,\phi_1)=\{(x,t)\in M \times ]0,T[\text{ such that }
u(x,t)=\breve u(x,t)\}.
$$ 
\item
At each point $(x,t)\in \mT( \phi_0,\phi_1)$, the functions $u$ and $\breve u$
are differentiable, and satisfy $du(x,t)=d\breve u(x,t)$.
In addition, the mapping $(x,t)\lmto du(x,t)$ is locally Lipschitz
on $\mT(\phi_0,\phi_1)$.
\item
If $\gamma(t)\in \mF(\phi_0,\phi_1)$, then
$\partial_xu(\gamma(t),t)=\partial_v L(\gamma(t),\dot \gamma(t),t)$.
As a consequence,
the set 
$$
\mT^*(\phi_0,\phi_1)
:=\{
(x,p,t)\in T^*M\times ]0,T[
\text{ such that }
(x,t)\in \mT
\text{ and } p=\partial_x u(x,t)= \partial_x\breve u(x,t)\}
$$
is invariant under the Hamiltonian flow,
and the restriction to $\tilde \mT(\phi_0,\phi_1)$
of the projection $\pi$ is a bi-locally-Lipschitz 
homeomorphism onto its image $\mT(\phi_0,\phi_1)$.
\end{enumerate}
\end{prop}
\proof
Let us fix a point $(x,t)\in M\times ]0,T[$.
There exist points $y$ and $z$ in $M$ such that
$u(x,t)=\phi_0(y)+c_0^t(y,x)$
and 
$\breve u(x,t)=\phi_1(z)-c_t^T(x,z)$,
so that 
$$
u(x,t)-\breve u(x,t)=\phi_0(y)-\phi_1(z)+c_0^t(y,x)+c_t^T(x,z)
$$
$$
\geq c_0^T(y,z)-(\phi_1(z)-\phi_0(y))
\geq 0.
$$
In case of equality, we must have 
$c_0^T(y,z)=c_0^t(y,x)+c_t^T(x,z)$.
Let 
$\gamma_1(s)\in C^2([0,t],M)$ satisfy
$\gamma_1(0)=y$, $\gamma_1(t)=x$ and
$\int_0^t L(\gamma_1(s),\dot \gamma_1(s),s)ds=c_0^t(y,x)$,
and 
let $\gamma_2(s)\in C^2([t,T],M)$ satisfy
$\gamma_2(t)=x$, $\gamma_2(T)=z$ and
$\int_0^tL(\gamma_2(s),\dot \gamma_2(s),s)ds=c_t^T(x,z)$.
The curve $\gamma:[0,T]\lto M$ obtained by pasting 
$\gamma_1$ and $\gamma_2$ clearly satisfies
$\int_0^TL(\gamma(s),\dot \gamma(s),s) ds=c_0^T(y,z)$,
it is thus a $C^2$ minimizer, and belongs to 
$\mF(\phi_0,\phi_1)$. As a consequence,
we have $(x,t)\in \mT(\phi_0,\phi_1)$. 

Conversely, we have:
\begin{lem}\label{subequ}
If $v$ is a viscosity sub-solution of $(HJ)$
satisfying $v_0=\phi_0$ and $v_T=\phi_1$,
then $\breve u \leq v \leq u$.
If $(x,t)\in \mT(\phi_0,\phi_1)$, then we have 
$v(x,t)=u(x,t)$.
\end{lem}

\proof
The inequality  $\breve u \leq v \leq u$ is easy.
For example, for a given point $(x,t)$
 there exists $y$ in $M$ such that
$u(x,t)=\phi_0(y)+c_0^t(y,x)$,
and for this value of $y$, we have 
$v(x,t)\leq \phi_0(y)+c_0^t(y,x)$,
hence $v(x,t)\leq u(x,t)$.
The proof that $\breve u\leq v$ is similar.
In order to prove the second part of the lemma,
it is enough to prove that 
$v(\gamma(t),t)=u(\gamma(t),t)$
for each curve $\gamma\in \mF(\phi_0,\phi_1)$.
Since $v$ is a sub-solution, we have 
$$
v(\gamma(T),T)\leq v(\gamma(t),t)+c_t^T(\gamma(t),\gamma(T)).
$$
On the other hand, we have
$$
 v(\gamma(t),t)\leq u(\gamma(t),t)
\leq u(\gamma(0),0)+c_0^t(\gamma(0),\gamma(t)).
$$
As a consequence of all these inequalities, we have 
$$
\phi_1(\gamma(T))=v(\gamma(T),T)
\leq u(\gamma(0),0)+c_0^t(\gamma(0),\gamma(t))+c_t^T(\gamma(t),\gamma(T))
\leq \phi_0(\gamma(0))+c_0^T(\gamma(0),\gamma(T))
$$
which is an equality because $\gamma \in \mF(\phi_0,\phi_1)$.
Hence all the inequalities involved are equalities, 
and we have 
$
v(\gamma(t),t)=u(\gamma(t),t).
$
\qed
The end of the proof of the proposition is straightforward.
Point 2 follows from Proposition \ref{Fathilip} applied
to the locally semi-concave functions $u$ and $-\breve u$.
Point 3 follows from Proposition \ref{proximal}.
\qed

\subsection{Optimal $C^1$ sub-solution}

The following result,
on which a large part of  the present paper is  based,
is inspired  from  \cite{FaSi}, but seems new in the present 
context.
%
%
%
%
%
%
\begin{prop}\label{soussolC1}
We have 
$$
C_0^T(\mu_0,\mu_T)
=\max _{v}
\Big(
\int_M v_Td\mu_T-\int_M v_0d\mu_0
\Big),
$$
where the maximum is taken on the set
of Lipschitz  functions $v:M\times [0,T]\lto \Rm$
which are $C^1$ on $M\times ]0,T[$ and 
satisfy 
the inequality
\begin{equation}\label{soussol}
\partial_tv(x,t)+H(x,\partial_x v(x,t),t)\leq 0
\end{equation}
at each point $(x,t)\in M\times ]0,T[.$
\end{prop}

\proof
First, let $v(x,t)$ be a continuous function of $M\times[0,T]$
which is differentiable on $M\times]0,T[$, where it satisfies 
(\ref{soussol}).
We then have, for each $C^1$ curve $\gamma(t):[0,T]\lto M$,
the inequality
$$
\int_0^T L(\gamma(t),\dot \gamma(t),t)dt
\geq
\int _0^T \partial_x v(\gamma(t),t)\cdot\dot \gamma(t)
-
H(\gamma(t), v(\gamma(t),t),t)dt
$$
$$
\geq
\int_0^T  \partial_x v(\gamma(t),t)\cdot\dot \gamma(t)
+\partial_t v(\gamma(t),t)dt
=v(\gamma(T),T)-v(\gamma(0),0).
$$
As a consequence, we get $v(y,T)-v(x,0)\leq c_0^T(x,y)$
for each $x$ and $y$,
so that 
$$
\int v_Td\mu_T-\int v_0d\mu_0\leq C_0^T(\mu_0,\mu_T).
$$
The converse follows directly from the next theorem,
which is an analog in our context of the main result of   \cite{FaSi}. 
\qed

\begin{thm}\label{construction of v}
For each Kantorovich optimal
pair $(\phi_0,\phi_1)$, 
there exists a  Lipschitz function 
$v:M\times[0,T]\lto \Rm$ which is $C^1$
on $M \times ]0,T[$, which  coincides 
with $u$ 
on $M\times\{0,T\}\cup  \mT(\phi_0,\phi_1)$,
and which satisfies the inequality (\ref{soussol}) strictly
at each point of $M \times ]0,T[-\mT(\phi_0,\phi_1)$.
\end{thm}
\proof
The proof of  \cite{FaSi} can't be translated to
our context in a straightforward way.
Our proof is different, and, we believe, simpler.
It is  based on:
\begin{prop}
There exists a function $V\in C^2(M\times[0,T],\Rm)$
which is null on $\mT(\phi_0,\phi_1)$ and which
is positive on $M\times ]0,T[-\mT(\phi_0,\phi_1)$,
and such that 
\begin{equation}\label{V}
\phi_1(y)=\min_{\gamma(T)=y}
\phi_0(\gamma(0))+
\int_0^T L(\gamma(t),\dot \gamma(t),t)-V(\gamma(t),t) dt.
\end{equation}
\end{prop}

\proof
Let us define the norm 
$$
\|u\|_2=\sum _{\theta\in \Theta}
 \|u\circ \theta\|_{C^2(B_1\times [0,T],\Rm)}
$$
of functions $u\in C^2(M\times[0,T],\Rm)$, where $\Theta$
is the atlas of $M$ defined in the Appendix.
Let us denote by $U$ the open set  $M\times ]0,T[-\mT(\phi_0,\phi_1)$.
We need a Lemma.

\begin{lem}
Let $U_1\subset U$ be an open set whose closure 
$\bar U_1$ is compact and contained in $U$, and let $\epsilon>0$
be given. There exists a function 
$V_1\in C^2(M\times[0,T],\Rm)$, 
which is positive on $U_1$ and null outside of $\bar U_1$
which is such that the equality (\ref{V}) holds with $V=V_1$,
and such that  $\|V_1\|_2\leq \epsilon$.
\end{lem}
\proof
Let us fix the open set $U_1$, the pair $(\phi_0,\phi_1)$ and $y\in M$.
We claim that the minimum in 
$$\min_{\gamma(T)=y}
\phi_0(\gamma(0))+
\int_0^T L(\gamma(t),\dot \gamma(t),t)-V_1(\gamma(t),t) dt$$
is reached at a path $\gamma$ the graph of which does not meet $U_1$, 
provided that the function $V_1$
is supported in $\overline U_1$ and is sufficiently small in
the $C^0$ topology.
In order to prove the claim, suppose the contrary.
There exists a sequence $V^1_n$ ($n\in\NN$) and a sequence $\gamma_n$ 
such that 
$$
\min_{\gamma(T)=y}
\phi_0(\gamma(0))+
\int_0^T L(\gamma(t),\dot \gamma(t),t)-V_n(\gamma(t),t) dt
$$
is reached at $\gamma_n$, the graph of $\gamma_n$ meets $U_1$,
$V_n$ is supported in $\overline U_1$ (for all $n\in \NN$)
and $V_n\rightarrow 0$ in the $C^0$ topology.
As a consequence each $\gamma_n$ is $C^2$ and the sequence 
$\gamma_n$ ($n\in\NN$) is a minimizing
sequence for 
\begin{equation}\label{eq: minimization}
\phi_1(y)=\min_{\gamma(T)=y}
\phi_0(\gamma(0))+
\int_0^T L(\gamma(t),\dot \gamma(t),t)dt.
\end{equation}
Hence this sequence is compact for the $C^2$ topology and, 
extracting a subsequence if needed, can be assumed to converge to some
$\gamma_{\infty}$. Clearly $\gamma_\infty$ is a minimizer for
\eqref{eq: minimization} the graph of which meets $ \overline U_1$.
This is a contradiction with $\overline U_1\subset U=M\times]0,T[
\backslash\mT(\phi_0,\phi_1)$
and the fact that the graph
of $\gamma_\infty$ is included in $\overline{\mT(\phi_0,\phi_1)}$
(see Definition \ref{dfn: T}).
 \qed

Let $U_n\subset U,n\in \Nm$ be a countable sequence of open sets
covering  $U$ and whose closures $\bar U_n$
are contained in $U$.
There exists a 
sequence $V_n$ of functions of  $C^2(M\times [0,T],\Rm)$
such that, for each $n\in \Nm$ :
\begin{itemize}
\item The function $V_n$ is positive in $U_n$ 
and null outside of $\bar U_n$. 
\item 
We have $\|V_n\|_2\leq 2^{-n}\epsilon$.
\item The equality (\ref{V}) holds for 
the function 
$V^n=\sum_{i=1}^n V_i$.
\end{itemize}
Such a sequence can be build inductively by applying the lemma
to the Lagrangian $L-V^{n-1}$ with 
$\epsilon_n=2^{-n}\epsilon$.
Since $\|V_n\|\leq 2^{-n}\epsilon$, the sequence $V^n$ is converging
in $C^2$ norm to a limit $V\in C^2(M\times [0,T],\Rm)$.
This function $V$ satisfies the desired properties.
The proposition is proved.
\qed

In order to finish the proof of the theorem, 
we shall consider the new Lagrangian $\tilde L=L-V$,
and the associated Hamiltonian $\tilde H=H+V$,
as well as the associated cost functions $\tilde c_s^t$.
Let 
$$
\tilde u(x,t):= \min_{y\in M}\phi_0(y)+\tilde c_0^t(y,x),
$$
be   the viscosity solution of the Hamilton-Jacobi equation
\begin{equation}\tag{$\tilde H J$}
\partial_t\tilde u+H(x,\partial_x \tilde u,t)=-V(x,t)
\end{equation}
emanating from $\phi_0$.
The equality (\ref{V}) says that $\tilde u_T=\phi_1=u_T$.
The function $\tilde u$ is 
Lipschitz  on $M\times[0,T]$, 
as a viscosity solution 
of $(\tilde HJ)$ emanating from a Lipschitz  function.
It is obviously a viscosity sub-solution of the equation 
$(HJ)$, which is strict outside of 
$M\times\{0,T\}\cup \mT(\phi_0,\phi_1)$
(where $V$ is positive).
This means that the inequality (\ref{soussol})
is strict at each point of differentiability of 
$\tilde u$  outside of 
$M\times\{0,T\}\cup \mT(\phi_0,\phi_1)$.
We have $\breve u\leq \tilde u\leq u$,
this relation being satisfied by each viscosity sub-solution
of $(HJ)$ which satisfies $u_0=\phi_0$ and $u_T=\phi_1$.
As a consequence, we have  $\breve u=\tilde u= u$ on $\mT(\phi_0,\phi_1)$,
and the function
$\tilde u$ is differentiable at each point 
of $\mT(\phi_0,\phi_1)$.
Furthermore,  we have $du=d\tilde u=d\breve u$ on this set.
We then obtain the desired function $v$ of the theorem 
from the function $\tilde u$ by regularisation,
applying  Theorem 9.2 of \cite{FaSi}. 
\qed

\section{Optimal objects of the direct problems}\label{ODP}
We prove Theorem A as well as the results of section \ref{problems}.
The following lemma generalizes a result of Benamou and Brenier,
see \cite{BeBr}.
\begin{lem}
We have the equality
$$
C_0^T(\mu_0,\mu_T)
= 
\min_{m_0\in \mI(\mu_0,\mu_T)} A(m_0)
=
\min_{m\in \mM(\mu_0,\mu_T)} A(m)
= 
\min_{\chi\in \mC(\mu_0,\mu_T)} A(\chi).
$$
Moreover $\chi(dv)=A(\chi)$ for every optimal $\chi$, where
$v$ is given by Theorem \ref{construction of v}
\end{lem}

\proof 
In view of Lemma \ref{inequalities}, it is enough to prove 
that, for each transport current $\chi\in \mC(\mu_0,\mu_T)$, we have  
$A(\chi)\geq C_0^T(\mu_0,\mu_T)$.
Let $v:M\times [0,T]\lto \Rm$
 be a Lipschitz  sub-solution of (HJ) which is $C^1$
on $M\times]0,T[$, and such that $(v_0,v_T)$ is a
Kantorovich  optimal pair.
For each current $\chi \in \mC(\mu_0,\mu_T)$,
we have $A(\chi)\geq \chi(dv)=C_0^T(\mu_0,\mu_T)$, which ends the proof.
\qed

From now on we choose and fix: 
\begin{itemize}
\item
An optimal Kantorovich pair $(\phi_0,\phi_1)$.
\item
A  Lipschitz sub-solution  $v:M\times [0,T]\lto \Rm$
of the Hamilton-Jacobi equation which 
satisfies $v_0=\phi_0$ and $v_T=\phi_1$ and 
which is $C^1$ on 
$M\times]0,T[$.
\item
A bounded  vectorfield 
$X(x,t):M\times]0,T[\lto TM$ which is locally Lipschitz
and satisfies 
\begin{equation}\label{refX}
X(x,t)=\partial_pH(x,\partial_xv(x,t),t)
~~\hbox{on}~~ \mT(\phi_0,\phi_1).
\end{equation}
\end{itemize}

\subsection{Characterization of optimal currents.}\label{OTC}
\begin{itshape}
Each optimal transport current $\chi$
can be written 
$$\chi=(X,1)\wedge \mu_{\chi}, 
$$
with a measure $\mu_{\chi}$ concentrated on  $\mT(\phi_0,\phi_1)$.
The current $\chi$ is then Lipschitz regular, 
so that there exists a transport interpolation  $\mu_t,t\in[0,T]$
 such that $\mu_{\chi}=\mu_t\otimes dt$ (see Appendix)
and such that 
$\mu_t=(\Psi_s^t)_{\sharp}\mu_s$ for each $s$ and $t$ in $]0,T[$.
\end{itshape}

\proof 
Let $\chi$ be an optimal transport current, that is a transport current
$\chi\in \mC(\mu_0,\mu_T)$ such that 
$A(\chi)=C_0^T(\mu_0,\mu_T)$.
Let us recall the definition of
the action  $A(\chi)$ that will be used here:
$$
A(\chi)=
\sup _{\omega\in \Omega^0} \Big(\chi(\omega^x,0)
-\int _{M\times[0,T]}
H(x,\omega^x(x,t),t)
d\mu_{\chi}\Big).
$$
Since $H(x,\partial_x v,t)+\partial_t v\leq 0$, we have 
$$
A(\chi)=\chi(dv)
\leq 
\chi(dv)
-\int 
H(x,\partial_x v(x,t),t)+\partial_t v
d\mu_{\chi}
=
\chi(\partial_x v,0)
-\int 
H(x,\partial_x v(x,t),t)
d\mu_{\chi}.$$
The other inequality holds by the  definition of $A$,  so that  
$$
\chi(dv)
=
\chi(dv)
-\int 
H(x,\partial_x v(x,t),t)+\partial_t v
d\mu_{\chi}=
\chi(\partial_x v,0)
-\int 
H(x,\partial_x v(x,t),t)
d\mu_{\chi},
$$
and we conclude that the function
$H(x,\partial_x v(x,t),t)+\partial_t v$
vanishes on the support of $\mu_{\chi}$,
or in other words that the measure $\mu_{\chi}$
is concentrated on the set $\mT(\phi_0,\phi_1)$.
In addition, for all form $\omega=\omega^x+\omega^tdt$,  
we have 
$$\chi(\partial_x v+\omega^x,0)-
\int H(x,\partial_x v +\omega^x,t)d\mu_{\chi}
\leq \chi(\partial_x v,0)-
\int H(x,\partial_x v,t)d\mu_{\chi}
=A(\chi).
$$
Hence the equality 
$$
\chi(\omega^x,0)=\int \partial_pH(x,\partial_x v ,t)(\omega ^x)d\mu_{\chi}
$$
holds for each form $\omega$.
This equality can be rewritten
$$
\chi(\omega)=\int \partial_pH(x,\partial_x v ,t)(\omega ^x)
+\omega^t 
d\mu_{\chi}
$$
which is precisely saying that 
$$
\chi=
(\partial_pH (x,\partial_x v(x,t),t),1)\wedge \mu_{\chi}
=
(X,1)\wedge \mu_{\chi}.
$$
The last equality follows from the fact that 
the vectorfields $X$ and $\partial_pH (x,\partial_x v(x,t),t)$
are equal on the support of $\mu_{\chi}$.
By the structure of Lipschitz regular transport currents,
we obtain the existence of a continuous family 
$\mu_t, t\in[0,T] $ of probability measures such that  
 $\mu_{\chi}=\mu_t\otimes dt$
and such that 
$\mu_t=(\Psi_s^t)_{\sharp}\mu_s$ for each $s$ and $t$ in $]0,T[$.
Since the restriction to a subinterval
$[s,t]\subset [0,T]$ of an optimal transport current
$\chi$
is clearly an optimal transport current for the transportation problem
between $\mu_s$ and $\mu_t$ with cost $c_s^t$, we obtain that 
the path $\mu_t$ is a transport interpolation.
\qed

\subsection{Characterization of transport interpolations.}
\begin{itshape}
Each transport interpolation $\mu_t$
satisfies 
$$\mu_t=(\Psi_s^t)_{\sharp}\mu_s$$
for each $(s,t)\in ]0,T[^2$.
The mapping
$$\mu_t \lmto (X,1)\wedge (\mu_t\otimes dt)
$$
is a bijection between the set of transport interpolations
and the set of optimal transport currents.
\end{itshape}
\vs

\proof
We fix a transport interpolation $\mu_t$ and two  times 
$s<s'$ in  $]0,T[$.
Let $\chi_1$ be a transport current 
on $M\times[0,s]$ between the measures $\mu_0$ and
$\mu_{s}$ which is optimal for the cost $c_0^{s}$,
let $\chi_2$ be a transport current
on $M\times [s,s']$ between the measures $\mu_{s}$
and $\mu_{s'}$ which is optimal for the cost $c_s^{s'}$
and let $\chi_3$ be a transport current
on $M\times [s',T]$ between the measures $\mu_{s'}$
and $\mu_T$ which is optimal for the cost $c_{s'}^{T}$.
Then the current $\chi $ on $M\times[0,T]$
which coincides with 
$\chi_1$ on $M\times[0,s]$, with 
$\chi_2$ on $M\times[s,s']$
and with $\chi_3$ on $[s',T]$ 
belongs to 
$\mC(\mu_0,\mu_T)$.
In addition, since $\mu_t$ is a transport interpolation, we have  
$$
A(\chi)=
C_0^{s}(\mu_0,\mu_{s})
+
C_{s}^{s'}(\mu_{s},\mu_{s'})
+
C_{s'}^T(\mu_{s'},\mu_T)
=C_0^T(\mu_0,\mu_T).
$$
Hence $\chi$ is an optimal transport current for the cost $c_0^T$.
In view of the characterisation of optimal
currents, this implies  that $\chi=(X,1)\wedge \mu_{\chi}$,
and that 
$$
\mu_{\chi}=\big( (\Psi_{s}^t)_{\sharp}\mu_{s}\big) \otimes dt
= \big((\Psi_{s'}^t)_{\sharp}\mu_{s'}\big) \otimes dt.
$$
By uniqueness of the continuous desintegration of 
$\mu_{\chi}$, we  obtain that, for each $t\in]0,T[$,
$
(\Psi_s^t)_{\sharp}{\mu_s}=
(\Psi_{s'}^t)_{\sharp}{\mu_{s'}},
$
and since this holds for all $s$ and $s'$,
that $(\Psi_s^t)_{\sharp}\mu_s=\mu_t$
for all $(s,t)\in ]0,T[^2$.
It follows that  $\chi=(X,1)\wedge (\mu_t\otimes dt)$.
We have proved that the mapping 
$$
\mu_t\lmto (X,1)\wedge (\mu_t\otimes dt)
$$
associates an optimal transport current to each transport interpolation.
This mapping is obviously injective, and it is surjective in view
of the characterization of optimal currents.
\qed

\subsection{Characterization of optimal measures.}\label{COM}
\begin{itshape}
The mapping 
$$\chi\lmto (X\times \tau)_{\sharp}\mu_{\chi}
$$
is a bijection between the set of optimal transport currents
and the set of optimal transport measures
($\tau:M\times[0,T]\rightarrow [0,T]$ is the projection on the second factor;
see Appendix).
Each optimal transport measure is thus invariant 
(see \eqref{X and flow} and Definition 
\ref{defn invariant}).
The mapping
$$m_0\lmto \mu_t=(\pi\circ \psi_0^t)_{\sharp}m_0
$$
is a bijection between the set of optimal initial measures
$m_0$ and the set of interpolations.
An invariant measure $m$ is optimal if and only if it is supported on
the set $\tilde \mT(\phi_0,\phi_1)$.\vs
\end{itshape}

\proof
If $m$ is an optimal transport measure, then
the associated current $\chi_m$ is an optimal transport
current, and $A(m)=A(\chi_m)$.
Let $\mu_m$ be the time component of $\chi_m$,
which is also the measure 
$(\pi\times\tau)_{\sharp}m$.
In view of
the characterization of optimal currents, 
we have $\chi_m=(X,1)\wedge \mu_m$.
We claim  that the equality 
$A(\chi_m)=A(m)$ implies that $m$ is supported
on the graph of $X$.
Indeed, we have 
the pointwise inequality 
\begin{equation}\label{legendre}
\partial_xv(x,t)\cdot V-H(x,\partial_xv(x,t),t)\leq L(x,V,t)
\end{equation}
for each $(x,V,t)\in TM\times]0,T[$.
Integrating with respect to $m$, we get the equality 
$$
A(\chi_m)=\chi_m(dv)
=\int_{TM\times[0,T]}
\partial_xv(x,t)\cdot V+\partial_tv(x,t)
dm(x,V,t)
$$
$$
=
\int _{TM\times[0,T]}
\partial_xv(x,t)\cdot V-H(x,\partial_xv(x,t),t)
dm(x,V,t)
= 
\int _{M\times[0,T]}
L(x,V,t)
dm(x,V,t)=
 A(m),
$$
which means that $m$ is concentrated
on the set where the inequality
(\ref{legendre})
is an equality, that is on the graph of
the vectorfield
$\partial_pH(x,\partial_xv(x,t),t)$.
Since $\mu_m$ is supported on $\mT$,
the measure $m$ is supported on $\tilde \mT$
and satisfies $m=(X\times \tau)_{\sharp}\mu_m$.
Let $\mu_t$ be the transport interpolation such that 
$\mu_m= \mu_t\otimes  dt$.
Setting $m_t=(X_t)_{\sharp}\mu_t$,
we have $m= m_t\otimes  dt$.
Observing that the relation
$$
X_t\circ \Psi_s^t =\psi_s^t\circ X_s
$$
holds on $\mT_s$, we conclude,
since $\mu_s$ is supported on $\mT_s$,
that 
$$
(\psi_s^t)_{\sharp} m_s=m_t,
$$
which means that the measure $m$ is invariant.

Conversely, let $m=m_t\otimes  dt$ be an invariant measure
supported on $\tilde \mT(\phi_0,\phi_1)$. 
We have 
$$
A(m)=\int_0^T \int_{TM} L (x,v,t)dm_t(x,v) dt
=
\int_0^T\int _{TM} L((\psi_0^t(x,v),t) dm_0(x,v) dt,
$$
and by Fubini,
$$
A(m)=\int_{TM} \int_0^T L((\psi_0^t(x,v),t)  dt dm_0(x,v)
=\int _{TM}\phi_1(\pi \circ \psi_0^T(x,v))-\phi_0(x) dm_0(x,v),
$$
and since $m_0$ is an initial transport measure, we get 
$$A(m)=\int_{TM} \phi_1 d\mu_T- \int _{TM}\phi_0 d\mu_0
=C_0^T(\mu_0,\mu_T).$$
\qed

%
%
%

%

\section{Absolute continuity}\label{AC}
In this section, we make the additional assumption that the initial
measure $\mu_0$ is absolutely continuous,
and  prove Theorem B.
The following lemma answers a question asked to us by Cedric Villani.

\begin{lem}
If $\mu_0$ or $\mu_T$ is absolutely continuous with
respect to the Lebesgue class,
then each interpolating measure  $\mu_t,t\in ]0,T[$, is absolutely continuous.
\end{lem}
\proof
If $\mu_t,t\in [0,T]$ is a transport interpolation,
we have proved that 
$$
\mu_t=(\pi\circ \psi_s^t\circ X_s)_{\sharp}\mu_s
$$
for each $s\in ]0,T[$, and $t\in [0,T]$.
Since the function $\pi\circ \psi_t^s\circ X_t$
is Lipschitz, it maps Lebesgue zero measure sets
into Lebesgue zero measure sets, and so 
it transport singular measures into singular measures.
It follows that if, for some $s\in ]0,T[$, 
the measure $\mu_s$ is not absolutely continuous,
then none of the measures $\mu_t,t\in[0,T]$ are absolutely
continuous. 
\qed

In order to continue the investigation of the specific properties
satisfied when $\mu_0$ is absolutely continuous, we first 
need some more general results.
Let $(\phi_0,\phi_1)$ be an optimal Kantorovich pair
for the measures $\mu_0$ and $\mu_T$ and for the cost 
$c_0^T$.
Recall that we have defined 
$\mF(\phi_0,\phi_1)\subset C^2([0,T],M)$ 
as the the set of curves $\gamma(t)$
such that 
$$
\phi_1(\gamma(T))=\phi_0(\gamma(0))+\int_0^TL(\gamma(t),\dot \gamma(t),t)dt.
$$ 
Let $\mF_0(\phi_0,\phi_1)$ be the set of initial velocities 
$(x,v)\in TM$ such that the curve $t\lmto \pi \circ \psi_0^t(x,v)$ 
belongs to $\mF(\phi_0,\phi_1)$. Note that there is a natural bijection
between $\mF_0(\phi_0,\phi_1)$ and $\mF(\phi_0,\phi_1)$.

\begin{lem}\label{mF}
The set $\mF_0(\phi_0,\phi_1)$ is compact. The maps $\pi$ and 
$\pi \circ \psi_0^T:\mF_0(\phi_0,\phi_1)\lto M$
are surjective.
If $x$ is a point of differentiability of $\phi_0$, then 
the set $\pi^{-1}(x)\cap\mF_0(\phi_0,\phi_1)$ contains one and only one point.
There exists a Borel measurable set $\Sigma\subset M$ of full measure,
whose points are points of differentiability of $\phi_0$, and 
such that the map 
$$
x\lmto S(x)=\pi^{-1}(x)\cap\mF_0(\phi_0,\phi_1)$$
is Borel measurable on $\Sigma$.
\end{lem}

\proof
The compactness of $\mF_0(\phi_0,\phi_1)$ follows from the fact,
already mentioned, that the set of minimizing extremals 
$\gamma:[0,T]\lto M$ is compact for the $C^2$- topology.

It is equivalent to say that the projection 
$\pi$ restricted to $\mF_0(\phi_0,\phi_1)$
is surjective, and to say that, for each point $x\in M$,
there exists a curve emanating from $x$ in $\mF(\phi_0,\phi_1)$.
In order to build such curves, recall that 
$$
\phi_0(x)= \max_{\gamma}
\phi_1(\gamma(T))-\int_0^TL(\gamma(t,\dot \gamma(t),t) dt
$$
where the maximum  is taken on the set of 
curves which satisfy $\gamma(0)=x$.
Any maximizing  curve is then a curve of 
$\mF(\phi_0,\phi_1)$ which satisfies $\gamma(0)=x$.
In order to prove that the map $\pi \circ \psi_0^T$ restricted 
to $\mF_0(\phi_0,\phi_1)$
is surjective, it is sufficient to build,
for each point $x\in M$, a curve in $\mF(\phi_0,\phi_1)$
which ends at $x$.
Such a curve is obtained as a minimizer in the expression
$$
\phi_1(x)=\min_{\gamma} \phi_0(\gamma(0))+
\int_0^TL(\gamma(t,\dot \gamma(t),t) dt.
$$

Now let us consider a point $x$ of differentiability of
$\phi_0$.
Applying the general result on the differentiability of viscosity solutions
to the Backward viscosity solution $\breve u$, we get that 
there exists a unique maximizer  to the problem 
$$
\phi_0(x)= \max_{\gamma}
\phi_1(\gamma(T))-\int_0^TL(\gamma(t,\dot \gamma(t),t) dt
$$
and that this maximizer is the extremal with 
initial condition 
$(x,\partial_pH(x,d\phi_0(x), 0))$.
As a consequence, there exists one and only one point $S(x)$
in $\mF_0(\phi_0,\phi_1)$ above $x$, and in addition we have the explicit expression
$$
S(x)=\partial_pH(x,d\phi_0(x), 0).
$$

Since the set of points of differentiability of $\phi_0$
has total Lebesgue measure --because $\phi_0$ is Lipschitz--
there exists a sequence $K_n$ of compact sets 
such that $\phi_0$ is differentiable at each point of $K_n$ 
 and such that the Lebesgue measure
of $M-K_n$ is converging to zero.
For each $n$, the set 
$\pi^{-1}(K_n)\cap \mF_0(\phi_0,\phi_1)$
is compact, and the restriction to this set 
of the canonical projection $\pi$ is injective and continuous.
It follows that the inverse function $S$ is continuous on $K_n$.
As a consequence, the map $S$ is Borel measurable on 
$\Sigma:= \cup_n K_n$.
\qed

\begin{lem}
The initial transport measure $m_0$ is optimal if and only if 
it is an initial  transport measure supported on $\mF_0(\phi_0,\phi_1)$.
\end{lem}

\proof
This statement is a reformulation of the result  in \ref{COM}
stating
that the optimal transport measures are the invariant measures
supported on $\tilde \mT(\phi_0,\phi_1)$.
\qed

\begin{prop}\label{ac}
If $\mu_0$ is absolutely continuous, then there exists 
a unique optimal initial measure $m_0$. 
There exists a Borel section 
$S:M\lto TM$ of the canonical
projection such that  $m_0=S_{\sharp}\mu_0$, this 
section is unique $\mu_0$-almost everywhere.
For each $t\in [0,T]$, the map
$\pi\circ \psi_0^t\circ S:M\lto M$
is then an optimal transport map between $\mu_0$ and $\mu_t$.
\end{prop}

\proof
Let $S:\Sigma \lto TM$ be the Borel map constructed in
Lemma  \ref{mF}.
For convenience,  we shall also denote by  $S$ 
the same map extended by zero outside 
of $\Sigma$, which is a Borel section
$S:M\lto TM$.
Since the set $\Sigma$ is of full Lebesgue measure,
and since the measure $\mu_0$ is absolutely continuous,
we have $\mu_0(\Sigma)=1$.
Let us consider the measure
$m_0=S_{\sharp}(\mu_{0|\Sigma})$.
This is a probability measure on 
$TM$, which is concentrated on 
$\mF_0(\phi_0,\phi_1)$, and which satisfies
$\pi_{\sharp}m_0=\mu_0$.
We claim that it is the only measure with these properties.
Indeed, if $\tilde m_0$ is a measure with 
these properties, then $\pi_{\sharp}\tilde m_0=\mu_0$,
hence the measure $\tilde m_0$ is concentrated
on $\pi^{-1}(\Sigma)\cap \mF_0(\phi_0,\phi_1)$.
But then, since 
$\pi$ induces a Borel isomorphism from 
$\pi^{-1}(\Sigma)\cap \mF_0(\phi_0,\phi_1)$ onto its image
$\Sigma$, of inverse $S$, we must have $\tilde m_0=S_{\sharp}\mu_0$.
As a consequence, the measure 
$m_0=S_{\sharp}\mu_0$ is the only candidate to be an 
optimal initial
transport measure.
Since we have already proved the existence of an
optimal initial  transport measure, it implies that 
$m_0$ is the only optimal initial transport measure.
Of course, we could prove directly that $m_0$ is
an initial transport measure, but as we have seen,
it is not necessary.
\qed

\subsection{Remark}
That there exists an optimal transport map 
if $\mu_0$ is
absolutely  continuous  could
be proved directly as a consequence of the following
properties of the cost function.
\begin{lem}\label{injectivity}
The cost function $c_0^T(x,y)$ is semi-concave on $M\times M$.
In addition, we have the following injectivity property for each $x\in M$:
If the differentials  $\partial_xc_0^T(x,y)$ and $\partial_xc_0^T(x,y')$
exist and are equal, then $y=y'$.
\end{lem}

In view of these properties of
the cost function, it is not hard to prove the following lemma
using a Kantorovich optimal pair in the 
spirit of works of Brenier \cite{Br:91}
and  Carlier \cite{Ca:un}.

\begin{lem}
There exists a compact subset $K\in M\times M$, 
such that the fiber $K_x=K\cap \pi_0^{-1}(x)$
contains one and only one  point for Lebesgue almost every $x$,
and which contains the support of all optimal plans.
\end{lem}
The proof of the existence of an optimal
map for an absolutely continuous measure $\mu_0$
can then be terminated using the following 
result, see \cite{Am:00}, Proposition 2.1.

\begin{prop}
A transport plan $\eta$ is induced from a transport map
if and only if it is concentrated on a $\eta$-measurable
graph.
\end{prop}

\subsection{Remark}
Assuming only that  $\mu_0$ vanishes on 
countably $(d-1)$-rectifiable sets,  we can conclude
that the same property holds for all interpolating measures 
$\mu_t, t<T$, and that Proposition \ref{ac} 
hold.
This is proved almost identically.
The only  refinement needed is that the set of singular points
of the semi-convex  function $\phi_0$ is a 
countably
$(d-1)$-rectifiable,
see \cite{CaSi}. 

%
%
%
\section{Aubry-Mather theory}\label{AMT}
We explain the relations between the results obtained so far
and Mather theory, and prove Theorem C.
Up to now, we have worked with fixed measures $\mu_0$ and $\mu_T$.
Let us study the optimal value $C_0^T(\mu_0,\mu_T)$ 
as a function of the measures $\mu_0$ and $\mu_T$.

\begin{lem}
The function 
$$(\mu_0,\mu_T)\lmto C_0^T(\mu_0,\mu_T)
$$
is convex and lower semi-continuous on the set of pairs of probability measures
on $M$.
\end{lem}
\proof
It follows directly from the expression
$$C_0^T(\mu_0,\mu_T)
=\max_{(\phi_0,\phi_1)}
\int_M \phi_1 d\mu_T-\int_M \phi_0 d\mu_0
$$
as a maximum of continuous linear functions.
\qed

From now on, we consider that the Lagrangian $L$ is defined 
for all times, $L\in C^2(TM\times\Rm,\Rm)$,
and satisfies
$$L(x,v,t+1)=L(x,v,t)$$
in addition to the standing hypotheses.
Let us restate  Theorem C with more details.
Recall that $\alpha$ is the action of Mather measures,
as defined in the introduction.\vs

\noindent
\textbf{Theorem C'.}
\begin{itshape}
There exists a Lipschitz vectorfield $X_0$ on $M$
such that all the Mather measures are supported on the graph of $X_0$.
We have 
$$
\alpha
=\min _{\mu} C_0^1(\mu,\mu),
$$
where the minimum is taken on the set of probability measures on $M$.
The mapping $m_0\lmto (\pi)_{\sharp}m_0$ is a bijection
between the set of Mather measures $m_0$ and the set
of probability measures $\mu$ on $M$ satisfying 
$C_0^1(\mu,\mu)=\alpha$.
More precisely, if $\mu$ is such a probability measure, 
then there exists one and only one initial transport measure
$m_0$ for the transport problem between $\mu_0=\mu$
and $\mu_1=\mu$ with cost $c_0^1$,
this measure is $m_0=(X_0)_{\sharp}\mu$, and 
it  is a Mather measure.
\end{itshape}\vs

The proof, and related digressions, occupy the end of the section.

\begin{lem}\label{alphaT}
The following minima
$$
\alpha_T:=\min_{\mu\in \mB_1(M)}\frac{1}{T}C_0^T(\mu,\mu), T\in \Nm
$$
exist and are all equal. In addition, any measure $\mu^1\in \mB_1(M)$
which is minimizing $C_0^1(\mu,\mu)$ is also minimizing 
$C_0^T(\mu,\mu)$ for all $T\in \Nm$.
\end{lem}

\proof
The existence of the minima follows from the compactness of
the set of probability measures and from the semi-continuity of
the function $C_0^T$.
Let $\mu^1$ be a minimizing measure for $\alpha_1$
and let $m^1$ be an optimal transport measure
for the transportation problem $C_0^1(\mu^1,\mu^1)$.
Let $m^T$ be the measure on $TM\times[0,T]$ obtained by concatenating
$T$ translated versions of $m^1$.
It means that $m^T$ is the only measure on $TM\times[0,T]$
whose restriction  to $TM\times[i,i+1]$
is obtained by translation from $m$, for each integer $i$.
It is easy to check that 
$m^T$ is indeed a transport measure between 
the $\mu_0=\mu^1$ and $\mu_T=\mu^1$ on the times interval
$[0,T]$, and that 
$A_0^T(m^T)=TA_0^1(m^1)$.
As a consequence, we have
$$
 T\alpha_T \leq C_0^T(\mu^1,\mu^1)
\leq A_0^T(m^T)
=TC_0^1(\mu^1,\mu^1)=T\alpha_1,
$$
which implies the inequality
$
\alpha_T\leq \alpha_1.
$

Let us now prove that $\alpha_T\geq \alpha_1$.
In order to do so, we consider an optimal measure $\mu^T$
for $\alpha_T$, 
and consider a transport  interpolation 
$\mu_t^T,t\in[0,T]$
between the measures $\mu_0=\mu^T$ and $\mu_T=\mu^T$.
Let us then consider, for $t\in[0,1]$, the measure
$$\tilde \mu^T_t:= \frac{1}{T}\sum_{i=0}^{T-1} \mu^T_{t+i},
$$
and note that 
$
T\tilde \mu^T_0=\mu^T_0+\sum_{i=1}^{T-1} \mu^T_{i}
=\mu^T_T+\sum_{i=1}^{T-1} \mu^T_{i}=T\tilde \mu^T_1.
$
In view of the convexity of the function $C_0^1$ 
$$
C_0^1(\tilde \mu^T_0,\tilde \mu^T_1)
=C_0^1\left(\frac 1 T \sum_{i=0}^{T-1}(\mu^T_{i},\mu^{T}_{i+1})\right)
\leq \frac{1}{T}\sum_{i=0}^{T-1} C_{i}^{i+1}(\mu_{i}^T,\mu_{i+1}^T)=
\frac{1}{T}C_0^T(\mu^T,\mu^T)=\alpha_T.
$$
Since $\tilde\mu^T_0=\tilde \mu^T_1$, 
this implies that $\alpha_1\leq \alpha_T$, as desired.
\qed

\begin{lem}
We have $\alpha_1\leq \alpha$.
\end{lem}

\proof
If $m_0$ is a Mather measure, 
then it is an initial measure for the transport
problem between $\mu_0=(\pi)_{\sharp}m_0$
and $\mu_1=(\pi)_{\sharp}m_0$
for the cost $c_0^1$.
As a consequence, we have 
$\alpha=A_0^1(m_0)\geq C_0^1(\mu_0,\mu_0)\geq \alpha_1.$
\qed

\begin{lem}\label{lem}
Let $\mu^1$ be a probability measure on $M$ 
such that $C_0^1(\mu^1,\mu^1)=\alpha_1$.
Then there exists a unique 
initial transport measure $m_0$ for the transportation problem
between $\mu_0=\mu^1$ and $\mu_1=\mu^1$ for the cost $c_0^1$.
This measure satisfies $(\psi_0^1)_{\sharp}m_0=m_0$.
We have $\alpha_1=A_0^1(m_0)\geq \alpha$, so that
$ \alpha=\alpha_1$
and
$m_0$ is a Mather measure.
There exists a constant $K$, which depends only on $L$,
such that the measure $m_0$ is supported on the 
graph  of a $K$-Lipschitz vectorfield.
\end{lem}

\proof
Let us fix a probability measure  $\mu^1$ 
on $M$ such that  $C_0^1(\mu^1,\mu^1)=\alpha_1$.
Let $X:M\times [0,2]\lto TM$
be a vectorfield associated to the transport problem 
$C_{0}^2(\mu^1,\mu^1)$ by Theorem A.
Note that $X_1$ is Lipschitz on $M$
with a Lipschitz constant $K$ which does not depend on $\mu_1$.
We choose $X$ once and for all and then fix it.

To each optimal transport measure  $m^1$
for the transport problem $C_0^1(\mu^1, \mu^1)$,
we associate the transport measure  $m^2$
on $TM\times[0,2]$
obtained by concatenation of two translated versions 
of $m^1$, as in the proof of Lemma \ref{alphaT}.
We have 
$$
A_0^2(m^2)=2A_0^1(m^1)=2\alpha_1=2\alpha_2=C_0^2(\mu^1,\mu^1).
$$
The measure $m^2$ is thus an optimal
transport measure for the transportation problem
$C_0^2(\mu^1,\mu^1)$.
Let $m_t,t\in [0,2]$ be the continuous family of probability measures
on $TM$ 
such that $m^2=m_t\otimes dt$. 
Note that 
$ m_t=(\psi_s^t)_{\sharp}  m_s$
for all $s$ and $t$ in $[0,2]$,
and that $m_0$ is the initial transport measure 
for the transportation problem $C_0^1(\mu^1,\mu^1)$ associated
to  $m^1$.
Since the measure $m^2$ was obtained by concatenation
of two translated versions of the same measure $m^1$, we must have  
$m_{t+1}= m_t$
for almost all $t\in]0,1[$, and,  by continuity,
 $m_0= m_1=m_2$. This implies that $m_0=(\psi_0^1)_{\sharp} m_0$.
Finally, the characterization of optimal measures implies
that 
$m_0=m_1=(X_1)_{\sharp}\mu^1$.
We have proved that the measure $(X_1)_{\sharp}\mu^1$
is the only optimal initial transport measure for the transportation
problem $C_0^1(\mu^1,\mu^1)$.
\qed

\noindent
\textsc{Proof of the theorem. }
Let $m_0$ be a Mather measure, and let $\mu_0=\pi_{\sharp}m_0$.
Note that we also have $\mu_0=(\pi\circ \psi_0^1)_{\sharp}m_0$.
As a consequence, $m_0$ is an initial transport measure for
the transport between $\mu_0$ and $\mu_0$ for the cost $c_0^1$,
and we have 
$$
\alpha=A_0^1(m_0)\geq C_0^1(\mu_0,\mu_0)\geq \alpha_1.
$$
Since $\alpha_1=\alpha$, all these inequalities are equalities,
so that $m_0$ is an optimal initial transport, and 
$C_0^1(\mu_0,\mu_0)= \alpha_1$.
It follows from Lemma \ref{lem} that $m_0$
is supported on the graph of a $K$-Lipschitz vectorfield.

Up to now, we have proved that each Mather measure is supported
on the graph of a $K$-Lipschitz vectorfield. There remains to prove that
all Mather measures are supported on a single $K$-Lipschitz graph.
In order to prove this, let us denote by $\tilde \mM\subset TM$ 
the union of the supports of Mather measures.
If $(x,v)$ and $(x',v')$ are two points of 
$\tilde \mM$, then there exists a Mather measure
$m_0$ whose support contains $(x,v)$
and a measure $m'_0$ whose support contains $(x',v')$.
But then the measure $(m_0+m'_0)/2$ is clearly a Mather measure 
whose support contains $\{(x,v),(x',v')\}$ and
is itself included in the graph of a $K$-Lipschitz vectorfield.
Assuming that $x$ and $x'$ lie in the image  $\theta(B_1)$
of a common chart, see appendix, so that $(x,v)=d\theta(X,V)$
and $(x',v')=d\theta(X',V')$, we obtain 
$$
\|V-V'\|\leq K\|x-x'\|.
$$
It follows that the restriction to $\tilde \mM$ of the canonical projection
$TM\lto M$ is a bi-Lipschitz homeomorphism, or equivalently
that the set $\tilde \mM$ is contained in the graph of a Lipschitz 
vectorfield. 
\qed

\appendix

\section{Notations and standing conventions}
\begin{itemize}
\item
$M$ is a compact manifold of dimension $d$, and $\pi:TM\lto M$
is the canonical projection.
\item
We denote by $\tau:TM\times[0,T]\lto [0,T]$
or $M\times[0,T]\lto [0,T]$ the projection on the second factor.
\item
If $N$ is any separable, complete, locally compact metric space 
(for example $M$,
$M\times[0,T]$,
$TM$ or $TM\times[0,T])$) the sets  
$\mB_1(N)\subset \mB_+(N) \subset \mB(N)$
are respectively the set of Borel probability measures,
non-negative Borel finite measures,
 and finite Borel signed measures.
If $C_c(N)$ is the set of continuous compactly supported  
functions on $N$, endowed with the topology of uniform convergence,
then the space $\mB(N)$ is identified with the set of continuous
linear forms on $C_c(N)$ by the Riesz theorem.
We will always endow the space  $\mB(N)$ with the weak-$*$ topology
that we will also call the weak topology.
Note that the set $\mB_1(N)$ is compact if $N$ is.
Prohorov's theorem states
that  a sequence of probability measures  $P_n\in \mB_1(N)$ 
has a  subsequence converging in $\mB_1(N)$ for the weak-$*$ topology  
if for all $\epsilon>0$ there exists a compact set
$K_\epsilon$ such that $P_n(N-K_\epsilon)\leq \epsilon$ for all $n\in\NN$.
See e.g. \cite{Vi:03,Du:02,Bi:99}.

\item
Given two manifolds $N$ and $N'$,  a Borel application 
$F:N\lto N'$, and a measure $\mu\in \mB(N)$, we define
the push-forward $F_{\sharp}\mu$ of $\mu$ by $F$ as the unique
measure on $N'$ which satisfies
$$
F_{\sharp}\mu(B)=\mu(F^{-1}(B))
$$
for all Borel set $B\in N$, or equivalently
$$
\int_{N'} fd(F_{\sharp}\mu)=\int_N f\circ F d\mu
$$
for all continuous function $f:N'\lto \Rm$.

\item
A family $\mu_t, t\in [0,T]$ of measures in $\mB(N)$ is called
measurable if the map
$t\lmto \int_N f_td\mu_t$ is Borel measurable
for each $f\in C_c(N\times [0,T])$.
We define the measure $\mu_t\otimes dt$ on $N\times [0,T]$
by 
$$
\int_{N\times [0,T]} fd(\mu_t\otimes dt)=
\int_0^T \int_N f_td\mu_t \, dt
$$
for each $f\in C_c(N\times[0,T])$.
The well-known desintegration theorem states that, if $\mu$
is a measure on $N\times [0,T]$ such that the projected measure
on $[0,T]$ is the Lebesgue measure $dt$, then there exists 
a measurable family of measures $\mu_t$ on $N$ such that 
$\mu=\mu_t\otimes dt$.

\item The set $\mK(\mu_0,\mu_T)$ of transport plans is defined
in section \ref{MKT}.
\item The set $\mI(\mu_0,\mu_T)$ of initial transport measures is defined
in section \ref{measures}.
\item The set $\mM(\mu_0,\mu_T)$ of transport measures is defined
in section \ref{measures}.
\item The set $\mC(\mu_0,\mu_T)$ of transport currents is defined
in section \ref{currents}.
\item
We  fix, once and for all, a finite atlas $\Theta$
of $M$, formed by charts $\theta:B_5\lto M$,
where $B_r$ is the open ball of radius $r$ 
centered at zero in $\Rm^d$. We  assume in addition
that the sets $\theta(B_1),\theta\in \Theta$ cover $M$.
\item
We  say that a vectorfield $X:M\lto TM$
is $K$-Lipschitz if, for each chart $\theta \in \Theta$,
the mapping  
$\Pi \circ(d\theta)^{-1}\circ X\circ \theta:B_5\lto  \Rm^d$
is $K$-Lipschitz on $B_1$, where $\Pi$ is the projection 
 $B_5\times\Rm^d\lto \Rm^d$.
\item
We mention the following results which are  used 
through the paper : 
There exists a  constant $C$ such that,
if $A$ is a subset of 
$M$, and $X_A:A\lto TM$ is a $K$-Lipschitz vectorfield,
then there exists a $CK$-Lipschitz vectorfield $X$ on $M$  which
extends $X_A$.
In addition, if $A$ is a subset of $M\times [0,T]$
 and $X_A:A\lto TM$ is a $K$-Lipschitz vectorfield,
then there exists a $CK$-Lipschitz vectorfield $X$ on $M\times[0,T]$  which
extends $X_A$.
If $A$ is a compact subset of $M\times [0,T]$
  and $X_A:A\cap M\times]0,T[\lto TM$ is a locally Lipschitz vectorfield
(which is $K(\epsilon)$-Lipschitz on $A\cap M\times [\epsilon,T-\epsilon]$),
then there exists a locally Lipschitz 
(which is $CK(\epsilon)$-Lipschitz on $M\times [\epsilon,T-\epsilon]$)
vectorfield $X$ on $M\times]0,T[$  which
extends $X_A$,
\end{itemize}

\small
\bibliographystyle{amsplain}
\providecommand{\bysame}{\leavevmode\hbox to3em{\hrulefill}\thinspace}

\end{document}